\documentclass[10.5pt,a4paper]{article}

\usepackage{amssymb}
\usepackage{amsmath}
\usepackage{tensor}
\usepackage{color}
\usepackage{stmaryrd}
 \usepackage{relsize}
\usepackage{cite}
\newtheorem{Th}{Theorem}[section]
\newtheorem{Prop}{Proposition}[section]
\newtheorem{Lm}{Lemma}[section]
\newtheorem{Co}{Corollary}

\newtheorem{Rm}{Remark}

\newcommand{\be}{\begin{equation}}
\newcommand{\ee}{\end{equation}}
\newcommand{\bes}{\begin{equation*}}
\newcommand{\ees}{\end{equation*}}
\newcommand{\R}{\mathbb{R}}

\newcommand\res{\mathop{\hbox{\vrule height 7pt width .5pt depth 0pt
\vrule height .5pt width 6pt depth 0pt}}\nolimits}

\def\theequation{\thesection.\arabic{equation}}
\def\theTh{\Roman{section}.\arabic{Th}}
\def\theProp{\Roman{section}.\arabic{Prop}}
\def\theCo{\Roman{section}.\arabic{Co}}
\def\theLm{\Roman{section}.\arabic{Lm}}

\def\theRm{\Roman{section}.\arabic{Rm}}
\newcommand{\reset}{\setcounter{equation}{0}\setcounter{Th}{0}\setcounter{Prop}{0}\setcounter{Co}{0}\setcounter{Lm}{0}\setcounter{Rm}{0}}

\def\al{\alpha}
\def\la{\lambda}

\def\pro{\pi_{\vec{n}}}
\def\bn{\vec{n}}

\def\bh{\vec{h}}
\def\ba{\vec{a}}

\def\bH{\vec{H}}
\def\bC{\vec{C}}
\def\bD{\vec{D}}
\def\bA{\vec{A}}
\def\bF{\vec{F}}
\def\bU{\vec{U}}

\def\bL{\vec{L}}

\def\bR{\vec{R}}

\def\bX{\vec{X}}
\def\bV{\vec{V}}

\def\bq{\vec{q}}

\def\bv{\vec{v}}
\def\bCC{\vec{C}}

\def\bB{\vec{B}}

\def\bp{\vec{\Phi}}

\def\bT{\vec{T}}

\def\bul{\bullet}

\def\res{\mathop{\hbox{\vrule height 7pt width .5pt 
depth 0pt\vrule height .5pt width 6pt depth 0pt}}\nolimits}

\begin{document}

\title{Structural Equations for Critical Points of\\ Conformally Invariant Curvature Energies in 4d }
\author{Yann Bernard\footnote{School of Mathematics, Monash University, 3800 Victoria, Australia.}}
\date{\today}
\maketitle

{\bf Abstract:} {\it This paper considers the Euler-Lagrange equations satisfied by the critical points of a large class of conformally invariant extrinsic energies for 4-manifolds immersed into $\R^{m\ge5}$. Using invariances and Noether's theorem, we convert the Euler-Lagrange equation in a system of equations with analytically favourable structures. The present paper generalises to the four-dimensional setting ideas originally developed by Tristan Rivi\`ere in his study of the Willmore energy in two dimensions.}\\



\reset

\section{Introduction and Main Results}


\subsection{Brief recap of the two-dimensional Willmore energy}\label{wil2d}

The Willmore energy is a fundamental conformally invariant functional in differential geometry, defined for smooth immersions of compact oriented surfaces into Euclidean space. For a surface \( \Sigma \subset \mathbb{R}^{m\ge3} \), its Willmore energy is given by
\[
\mathcal{W}(\Sigma) = \int_\Sigma |\bH|^2 \, dA,
\]
where \( \bH \) is the mean curvature vector and \( dA \) the induced area measure. This energy captures the total bending of the surface and remains invariant under Möbius transformations of \( \mathbb{R}^3 \cup \{\infty\} \). The conformal invariance of the Willmore energy becomes manifest when written
$$
\int_\Sigma|\bH|^2dA\;=\;\dfrac{1}{2}\int_\Sigma|\bh_0|^2dA+\int_\Sigma KdA\:,
$$
where $\bh_0$ denotes the trace-free part of the second fundamental form, while $K$ is the Gauss curvature. 
The first integral on the right is clearly conformally invariant, while the second integral is a topological constant by the Gauss-Bonnet theorem. \\

The historical origins of the Willmore energy are multifaceted: first introduced by Sophie Germain in the context of elasticity theory, where it models the bending energy of thin shells, it was studied by conformal geometer Blaschke and his students in an effort to combine minimal surfaces and conformal invariance (cf. \cite{Ber1} and the references therein). It was popularised by Thomas Willmore \cite{Wil} in the 1960s, and has since revealed itself as a key quantity appearing in numerous research papers coming from diverse areas of science and technology (see \cite{Ber1}). Summarising such a vast topic in a few paragraphs is necessarily doomed to bias. The next section is no exception.\\

\noindent 
One foundational result due to Li and Yau \cite{LY} establishes that for any compact immersed surface \( \Sigma \subset \mathbb{R}^{m} \), the Willmore energy satisfies the lower bound
\[
\mathcal{W}(\Sigma) \geq 4\pi,
\]
with equality if and only if \( \Sigma \) is a round sphere. When the energy is strictly below $8\pi$, the immersion must be an embedding. These statements confirm the rigidity of minimal bending configurations and serves as the starting point for much of the modern analysis of the Willmore functional. The celebrated Willmore conjecture, posed in the 1960s and resolved in 2014 by Marques and Neves \cite{MN}, asserts that any torus immersed in \( \mathbb{R}^3 \) must satisfy \( \mathcal{W}(\Sigma) \geq 2\pi^2 \), with equality uniquely for the stereographic image of the Clifford torus in \( S^3 \). Their resolution of the conjecture, using min-max theory in the Almgren--Pitts framework, marked a major milestone in geometric variational problems and introduced powerful new techniques for studying conformally invariant energies.\\

Beyond its variational significance, the Willmore functional serves as a bridge between extrinsic surface theory and deep structures in conformal geometry. In the setting of asymptotically hyperbolic manifolds, particularly those with Poincaré--Einstein metrics, the Willmore energy appears as a renormalised area of minimal surfaces asymptotic to a given boundary. This perspective, motivated by the AdS/CFT correspondence in theoretical physics, was developed by Graham and Witten \cite{GW}, and reveals that the Willmore energy of a boundary surface \( \Sigma \subset S^3 \) can be interpreted as the finite part in the asymptotic expansion of the area of a minimal surface in \( \mathbb{H}^4 \) with boundary \( \Sigma \). This realisation connects the Willmore functional to the renormalised volume of Einstein manifolds and embeds it into the broader framework of conformal invariants arising from geometric scattering theory and holographic correspondences.\\


The critical points of the Willmore energy -- be they minimisers or not -- are of course interesting in their own right. They are smooth immersions $\bp:\Sigma\rightarrow\mathbb{R}^m$ satisfying a quasilinear elliptic equation of order four (for the immersion):
\be\label{will2}
\Delta_\perp\bH+(\bH\cdot\bh^{ij})\bh_{ij}-2|\bH|^2\bH\;=\;\vec{0}\:.
\ee
Here $\bh$ is the second fundamental form, $\Delta_\perp:=\pro\nabla_j\pro\nabla^j$ is the covariant Laplacian in the normal bundle, and $\pro$ denotes projection onto the normal space. Repeated indices are summed over using the induced metric $g_{ij}=\nabla_i\bp\cdot\nabla_j\bp$. We note that minimal surfaces are examples.\\

\noindent
Obtaining information about solutions to (\ref{will2}) is a particularly challenging task as compared to second-order elliptic problems more commonly studied in differential geometry. In codimension 1, Kuwert and Sch\"atzle (see {\it inter alia} \cite{KS1, KS2}) used an ``ambient approach" developed by Simon \cite{Sim}. A different, ``parametric" approach suitable in all codimensions was devised by Rivi\`ere in \cite{Riv1}. That approach, relying on the existence of certain conservation laws, was later found \cite{Ber1} to stem from Noether's theorem. The idea consists in noting that (\ref{will2}) results from varying a translation-invariant energy and must therefore be in divergence form. In fact, one verifies directly using the Codazzi identity $\pro\nabla_i\bh^{ij}=2\pro\nabla^j\bH$:
\be\label{will22}
\nabla_i\big(\pro\nabla^i\bH+(\bH\cdot\bh^{ij}-|\bH|^2g^{ij})\nabla_j\bp\big)\;=\;\Delta_\perp\bH+(\bH\cdot\bh^{ij})\bh_{ij}-2|\bH|^2\bH\:.
\ee
For a critical point, an application of the Poincar\'e lemma yields the existence of a 1-vector-valued 2-form $\bL$ satisfying\footnote{The codifferential operator is understood with respect to the induced metric.}
\bes
(d^\star\bL)^i\;=\;(\bH\cdot\bh^{ij}-|\bH|^2g^{ij})\nabla_j\bp+\pro\nabla^i\bH\:.
\ees
This identity can be contracted on $d\bp$ in two different ways: using the scalar and wedge products in ambient space. Doing so yields
\be\label{will222}
(d^\star\bL)^i\cdot\nabla_i\bp\;=\;0\qquad\text{and}\qquad (d^\star\bL)^i\wedge\nabla_i\bp\;=\;\pro\nabla^i\bH\wedge\nabla_i\bp\:.
\ee
These equations were derived in \cite{Riv1}\footnote{In \cite{Riv1}, the 2-form $\bL$ is replaced by the real-valued function $\star\bL$.}. \\

An interesting question is whether (\ref{will222}) is equivalent to (\ref{will2}). The answer is ``not quite" and is explained in \cite{BR1}. When (\ref{will222}) holds for some 2-form $\bL$, then there exists a divergence-free, trace-free symmetric 2-tensor $T$ such that
\be\label{willPS}
\Delta_\perp\bH+(\bH\cdot\bh^{ij})\bh_{ij}-2|\bH|^2\bH\;=\;-T^{ij}\bh_{ij}\:.
\ee
With $T=0$, we recover the Willmore equation. It is shown in \cite{BR1} and \cite{BPP} that parallel mean curvature surfaces, i.e. those satisfying $\pro\nabla^j\bH=\vec{0}$ for all $j$, are solutions of (\ref{willPS}) with $T^{ij}=2\bH\cdot\bh_0^{ij}$. In this expression, $\bh_0:=\bh-\bH g$ is the trace-free second fundamental form. Per Codazzi, $T$ is indeed divergence-free when $\bH$ is parallel. \\

Although (\ref{will222}) is not quite equivalent to (\ref{will2}), it can nonetheless yield fruitful information. Note that
$$
(d^\star\bL)^i\cdot\nabla_i\bp\;=\;\nabla_i\big(\bL^{ij}\cdot\nabla_j\bp\big)\:,
$$
since $\bL^{ij}\cdot\bh_{ij}=0$. By the same token -- and a few elementary manipulations -- we can rewrite the second equation in (\ref{will222}). Altogether, we obtain
\bes
\nabla_i\big(\bL^{ij}\cdot\nabla_j\bp\big)\;=\;0\qquad\text{and}\qquad \nabla_i\big(\bL^{ij}\wedge\nabla_j\bp-\bH\wedge\nabla^i\bp\big)\;=\;\vec{0}\:.
\ees
Applying again the Poincar\'e lemma, we deduce the existence of a scalar-valued 2-form $S$ and a 2-vector-valued 2-form $\bR$ satisfying
\bes
(d^\star S)^i\;=\;\bL^{ij}\cdot\nabla_j\bp\qquad\text{and}\qquad (d^\star\bR)^i\;=\;\bL^{ij}\wedge\nabla_j\bp-\bH\wedge\nabla^i\bp\:.
\ees
These expressions can be recast via defining the 0-forms $\vec{\ell}=\star\bL$, $U=\star S$, and $\bV=\star\bR$. One finds
\bes
\nabla U\;=\;\vec{\ell}\cdot\nabla \bp\qquad\text{and}\qquad \nabla \bV\;=\;\vec{\ell}\wedge\nabla \bp-\bH\wedge\nabla^\perp\bp\:,
\ees
where $\nabla^\perp$ is the gradient rotated by 90 degrees. 
For algebraic reasons, these two equations are in fact linked (see \cite{Marq} for further elaboration on this point). Let $\vec{\eta}^{ij}:=\nabla^i\bp\wedge\nabla^j\bp$. One verifies that\footnote{The $\bullet$ notation is clarified in the Appendix. It is an operation which returns a 2-vector from a pair of 2-vectors.}
\begin{equation}\label{sysUV1}
\left\{\begin{array}{rcl}
\nabla^iU&=&\vec{\eta}^{\,ij}\cdot\nabla_j\bV\\[1ex]
-\nabla^i\bV&=&\vec{\eta}^{\,ij}\bullet\nabla_j\bV+\vec{\eta}^{\,ij}\nabla_jU\:.
\end{array}\right. 
\end{equation}
Of course, in dimension 2, we may always without loss of generality suppose that locally our induced metric is conformal: $g_{ij}=e^{2\la}\delta$, where $\delta$ denotes the flat metric. In this setting, the Gauss map $\bn$ satisfies $\star\bn=e^{-2\la}\partial_x\bp\wedge\partial_y\bp$. The 2-vector $\vec{\eta}$ is proportional to $\star\bn$. In this notation, the system (\ref{sysUV1}) takes the form given in Rivi\`ere's original work:
\begin{equation}\label{sysUV2}
\left\{\begin{array}{rcl}
\nabla^\perp U&=&(\star\bn)\cdot\nabla\bV\\[1ex]
-\nabla^\perp \bV&=&(\star\bn)\bullet\nabla\bV+(\star\bn)\nabla U\:.
\end{array}\right. 
\end{equation}
This remarkably symmetric and simple structure originally discovered by Rivi\`ere in \cite{Riv1} is the key that unlocks many analytical results and it has been exploited in numerous subsequent works, see {\it inter alia.} \cite{BLM, LR, Riv2}. However analytically amenable (\ref{sysUV2}) might be, there remains nonetheless the question of linking the potentials $U$ and $\bV$ back to geometric quantities. This was also done in \cite{Riv1}, namely:
$$
2\Delta\bp\;=\;\text{div}( U\cdot\nabla\bp+\bV\bullet\nabla\bp)\:.
$$
According to this narrative, in order to obtain analytical information for the mean curvature $\bH$, one needs analytical information for the potentials $U$ and $\bV$. Such information is at reach upon differentiating (\ref{sysUV1}). One uncovers a system of second-order elliptic PDEs with a Jacobian structure on the right-hand side. The good analytical dispositions of such a structure may be subsequently put to work, and with the help of this strategy, Rivi\`ere obtained energy estimates and several corollaries in \cite{Riv1}. The conservation laws for $\bL$, $S$, and $\bR$ can be recovered by evaluating the Noether fields associated with the translation, dilation, and rotation invariances of the Willmore energy \cite{Ber1}. Noether's method has further proved its worth when used in numerous related contexts, see {\it inter alia} \cite{BWW, LiY}.

\subsection{Conformally invariant energies in dimension four}\label{confinv}

In dimension 2, the Riemann tensor is literally reduced to the shadow of itself: its double contraction $2K$. In dimension 4, the situation dramatically complicates as the Riemann tensor now comes into full force and can be neither recovered from its double trace (as in 2d), nor from its trace (as in 3d). There is still, of course, in dimension 4 a topological invariant given by the Chern-Gauss-Bonnet theorem. Namely, letting $\text{Rm}$, $\text{Ric}$, and $R$ denote respectively the Riemann, Ricci, and scalar curvature tensors, the quantity
$$
\int_\Sigma\big(|\text{Rm}|^2-4|\text{Ric}|^2+R^2\big)d\text{vol}_g\;=\;32\pi^2\chi(M)
$$
is a topological invariant, with $\chi(M)$ denoting the Euler characteristic. But there is yet another intrinsic quantity of interest: the Weyl tensor. It is the totally traceless part of the Riemann tensor and is given by\footnote{$\owedge$ stands for the Kulkarni-Nomizu product.}
$$
W\;=\;\text{Rm}-\dfrac{1}{2}P\owedge g\qquad\text{where}\quad P\;:=\;\text{Ric}-\dfrac{R}{6}g
$$
is the Schouten tensor. The $(3,1)$-Weyl tensor has the property of being conformally invariant, and a celebrated result of Weyl states that a Riemannian manifold is locally conformally flat if and only if its Weyl tensor vanishes identically. Naturally, 
$$
\int_\Sigma|W|^2d\text{vol}_g
$$
constitutes a conformally invariant intrinsic energy. \\

In the submanifold setting, the diversity of conformally invariant energies increases, if only due to combinatorics. Indeed, from the pointwise conformal invariant tensor $\bh_0$ one can construct four distinct conformally invariant energies with integrands 
$$
|\bh_0|^4\quad,\quad \langle\bh_0^4\rangle:=\big(\bh_0^{ij}\cdot\bh_0^{kl}\big)\big((\bh_0)_{ij}\cdot(\bh_0)_{kl}\big)
$$
and
$$
|h_0^2|^2:=\big((\bh_0)^{ik}\cdot(\bh_0)_k^j\big)\big((\bh_0)_{il}\cdot(\bh_0)_j^l\big)\quad,\quad \text{Tr}(\bh_0^4):=\big((\bh_0)^{ij}\cdot(\bh_0)^{kl}\big)\big((\bh_0)_{il}\cdot(\bh_0)_{kj}\big)\:.
$$
The extrinsic expression of the Riemann tensor is well-known:
$$
\text{Rm}=\dfrac{1}{2}\bh\stackrel{\cdot}{\owedge}\bh\qquad\text{i.e.}\qquad R^{ijkl}\;=\;\bh^{ik}\cdot\bh^{jl}-\bh^{il}\cdot\bh^{jk}\:.
$$
This formulation yields that of the Weyl tensor in extrinsic terms:
$$
W\;=\;\dfrac{1}{2}\bh_0\stackrel{\cdot}{\owedge}\bh_0-\dfrac{1}{2}(h_0^2\owedge g)-\dfrac{1}{6}|\bh_0|^2g\owedge g\:,
$$
where, for notational convenience $(h_0^2)^{ij}:=\bh_0^{ik}\cdot(\bh_0)_k^j$. From this expression, one deduces easily that
$$
|W|^2\;=\;2\langle\bh_0^4\rangle-2\text{Tr}_g\bh_0^4-2|h_0^2|^2+\dfrac{1}{3}|\bh_0|^4\:.
$$
This confirms that $|W|^2$ is a linear combination of the four conformally invariant building blocks identified above. \\

In codimension 1, the four energies coincide pairwise. In an effort to construct more, one could adopt the approach suggested in \cite{MoN} and introduce such quantities as $\big(\text{Tr}(h_0^3)\big)^{4/3}$. Although indeed conformally invariant, these quantities can only be defined in codimension 1. That is too strong a restriction, and one can instead favour a different approach, which consists in noting that 
$$
\int_\Sigma|\pro d\bH|^2d\text{vol}_g
$$
scales correctly in dimension 4. It is not however a conformally invariant energy, but it can be modified to become one. Namely, the energy
$$
\mathcal{E}_A\;:=\;\int_\Sigma\big(|\pro d\bH|^2-|\bH\cdot\bh|^2+7|\bH|^4\big) d\text{vol}_g
$$
turns out to be conformally invariant! In codimension 1, this energy was first identified by Guven \cite{Guv}. The author derives it from scratch by evaluating and cleverly recasting the defect of conformal invariance of the Dirichlet energy of mean curvature. Following different routes, other authors have arrived at the energy $\mathcal{E}_A$. Graham and Reichert \cite{GR} on one hand, and Zhang \cite{Zha} on the other, confirmed the validity of $\mathcal{E}_A$ in all codimensions, using the correspondence between renormalised volume in Poincar\'e-Einstein metric and conformal invariants. Blitz, Gover, and Waldron \cite{BGW} use instead tractor calculus and produce numerous conformal invariant energies, in particular $\mathcal{E}_A$. In some references, $\mathcal{E}_A$ is rendered as $\mathcal{E}_{GR}$, for Graham-Reichert. \\

The Euler-Lagrange equation for $\mathcal{E}_A$ is quite a handful. It is derived in Section~\ref{vardh} (see also \cite{GR}):
\begin{eqnarray}\label{horrib}
&&\vec{\mathcal{W}}\;\;=\;\;\dfrac{1}{2}\Delta_\perp^2\bH+\dfrac{1}{2}\Delta_\perp\langle\bH\cdot\bh,\bh\rangle-7\Delta_\perp(|\bH|^2\bH)+8\pro\nabla_j\big(\bH\nabla^j|\bH|^2\big)+4\pro\nabla_j\big((\bH\cdot\bh^j_i)\pro\nabla^i\bH\big)\nonumber \\
&&\hspace{1.5cm}+\:\Big[\dfrac{1}{2}\bh^{ij}\cdot\Delta_\perp\bH-2\nabla_i\bH\cdot\pro\nabla_j\bH+2(\bH\cdot\bh^{ik})(\bH\cdot\bh^j_k)  \nonumber\\
&&\hspace{2cm}+\:\dfrac{1}{2}\langle\bH\cdot\bh,\bh\rangle\cdot \bh^{ij}-7|\bH|^2\bH\cdot\bh^{ij} +\big(|\pro d\bH|^2-|\bH\cdot\bh|^2+7|\bH|^4\big)g^{ij}\Big]\bh_{ij}\nonumber\\
&&\hspace{.5cm}\;\;=\;\;\vec{0}\:.
\end{eqnarray}
The energy $\mathcal{E}_A$ is, for many of its aspects, the analogue of the two-dimensional Willmore energy. Both are conformally invariant, both include minimal immersions as critical points, both have a quasilinear elliptic system (of order 4 in the immersion for 2d Willmore, and of order 6 in the four-dimensional setting). However, unlike the Willmore energy, $\mathcal{E}_A$ is not bounded from below. Recent works have shown that minimisers of $\mathcal{E}_A$ are hard to come by. Martino \cite{Mart} proves that $\mathcal{E}_A$ is unbounded from below on closed hypersurfaces. In an effort to thwart this annoying feature, one may add some other lower-order conformally invariant energy. For example, as is elementary to check \cite{Vya}, the energy $\mathcal{E}_A+\beta\int|\bh_0|^4$ is bounded below by $8\pi^2$ for $\beta>\frac{1}{12}$. In the recent work \cite{WY}, the authors establish a connected sum energy reduction for the energy $\mathcal{E}_A$.\\

It takes little imagination to accept that $\mathcal{E}_A$ is not the only conformally invariant energy of its kind and that there might at least be another energy involving this time 
$$
\int_\Sigma|\pro d\bh|^2d\text{vol}_g\:.
$$
An analogue of the following statement appears in \cite{BGW} as well as in \cite{Mart}. 
\begin{Prop}\label{EC}
The energy
$$
\mathcal{E}_C\;:=\;\int_\Sigma\big(|\pro d\bh|^2-12|\bH\cdot\bh|^2+6|\bH|^2|\bh|^2+60|\bH|^4\big)d\text{vol}_g
$$
is conformally invariant. 
\end{Prop}
$\hfill\square$\\

\noindent
As we will see in Section~\ref{UC}, the proof of this statement is a simple application of Simons' identity. In particular, we will obtain that
\be\label{simid}
16\mathcal{E}_A-\mathcal{E}_C+\int_\Sigma\big(3|h_0^2|^2-\langle\bh_0^4\rangle\big)d\text{vol}_g\;=\;16\pi^2\chi(\Sigma)-\dfrac{3}{2}\int_\Sigma|W|^2d\text{vol}_g\:.
\ee
The left-hand side appears extrinsic, but is in reality intrinsic. \\

In this paper, we will be interested in linear combinations of $\mathcal{E}_A$, $\mathcal{E}_C$, and any linear combination of the four conformally invariant energies of order $\mathcal{O}(|\bh_0|^4)$, which we generically denote by $\mathcal{E}_0$. We of course expound any topological invariant from the energy so as not to be distracted by their ghost Noether fields (which are complicated looking zeros). Owing to (\ref{simid}), we will thus be dealing with energies of so-called ``$\mathcal{E}_A$-type":
\be\label{typen}
\mathcal{E}\;=\;\mathcal{E}_A+\mathcal{E}_0\:.
\ee
From the point of view of Noether's fields and underlying structures, we will see that critical points of energies of the type (\ref{typen}) satisfy analogous equations to those found in the 2d Willmore problem. In particular, they satisfy a system of quasilinear elliptic equations of order six in the immersion. {\it A contrario}, energies of the type $\mathcal{E}_0$ satisfy a fully nonlinear system of equations whose study demands completely different methods. Of course, the most important representative of the $\mathcal{E}_0$ family is the intrinsic Bach energy $\int|W|^2$. Its critical points -- when varying the metric -- satisfy the Bach equation: 
$$
B^{ab}\;:=\;\Delta P^{ab}-\nabla_c\nabla^aP^{cb}+P_{ij}W^{aibj}\;=\;0\:.
$$
Exact solutions of the Bach equation include locally conformally flat spaces (for which $W=0$), Einstein manifolds, as well as self-dual and anti-self-dual manifolds (with $W=\pm\star W$). Analytical properties of Bach-flat manifolds have been considered in the literature, e.g. in \cite{TV, LM}. They will not be discussed in the present work. We content ourselves with recording a basic fact for our future use: the Bach tensor is a TT-tensor; it is divergence-free and trace-free. 

\subsection{Equation structures of $\mathcal{E}_A$-type energies}

We study in this section energies of the type $\mathcal{E}=\mathcal{E}_A+\mathcal{E}_0$, where $\mathcal{E}_0$ is any linear combination of the four conformally invariant energies of order $\mathcal{O}(|\bh_0|^4)$. We will closely follow the ideas summarised in section~\ref{wil2d}. As the Euler-Lagrange equation is quite unmanageable as written in (\ref{horrib}), the first thing to do is to express it in divergence-form. This can be accomplished by evaluating the Noether field corresponding to translation invariance. More precisely, in Section~\ref{Noether}, we show
\begin{Th}\label{cons}
A critical point of $\mathcal{E}$ satisfies the conservation law
$$
d^\star\bV\;=\;\vec{\mathcal{W}}\;=\;\vec{0}\:.
$$
The vector-field $\bV$ is given by
\bes
\bV^j\;=\;G^{ij}\nabla_i\bp-\pro\nabla_i\bF^{ij}+\bC^j\:.
\ees
In this expression\footnote{$\pro$ denotes projection on the normal subspace.}, we have set
\bes
\left\{\begin{array}{lcl}
G^{ij}&=&\dfrac{1}{2}\bh^{ij}\cdot\Delta_\perp\bH+ |\pro\nabla\bH|^2g^{ij}-2\nabla^j\bH\cdot\pro\nabla^i\bH+\mathcal{O}(|\bh|^4)\\[1ex]
\bF^{ij}&=&-\dfrac{1}{2}\Delta_\perp\bH\,g^{ij}+\mathcal{O}(|\bh|^3) \\[1ex]
\bC^j&=&-2(\bh^{ij}\cdot\nabla_i\bH)\bH+2(\bH\cdot\bh^{ij})\pro\nabla_i\bH\:.\end{array}\right.
\ees
Equivalently, 
\bes
\bV^j\;=\;\dfrac{1}{2}\nabla^j\Delta_\perp\bH+\nabla_kA^{kj}+\mathcal{O}\big(|\bh|^2|\pro d\bH|+|\pro d\bH|^2+|\bh|^4\big)\:,
\ees
with
$$
A^{kj}\;:=\;\mathcal{O}\big(|\bh||\pro d\bH|+|\bh|^3\big)\:.
$$
\end{Th}
$\hfill\square$\\

\noindent
As expected, the leading-order term in the Euler-Lagrange equation is linear and proportional to  $\Delta_\perp^2\bH$. \\

Integrating the equation $d^\star\bV=\vec{0}$ once yields the existence of a 1-vector-valued 2-form $\bL$ satisfying $d^\star\bL=\bV$. Contracting $d^\star\bL$ on $d\bp$ yields a pair of equations, henceforth referred to as the ``$\bL$-system":

\begin{Th}\label{ThsysL}
If $\bp$ is a critical point of the energy $\mathcal{E}$, there exists a 1-vector-valued 2-form $\bL$ satisfying
\begin{equation}\label{sysL20}
\left\{\begin{array}{lcl}(d^\star\bL)^j\cdot\nabla_j\bp&=&G^i_i\;\;=\;\;\Delta|\bH|^2\\[1ex]
(d^\star\bL)^{j}\wedge\nabla_j\bp&=&-\pro\nabla_i\bF^{ij}\wedge\nabla_j\bp+\bC^j\wedge\nabla_j\bp\:.
\end{array}\right.
\end{equation}
\end{Th}
$\hfill\square$\\

\noindent
This system is not quite equivalent to the original Euler-Lagrange equation, as the following statement outlines. 
\begin{Th}\label{ThPS}
There exists a 1-vector-valued 2-form $\bL$ satisfying (\ref{sysL20}) if and only if there exists a divergence-free, trace-free, symmetric 2-tensor $T$ such that
\be\label{bawo}
\vec{\mathcal{W}}\;=\;-\langle T,\bh\rangle\:.
\ee
\end{Th}
$\hfill\square$\\

\noindent
The case $T=0$ corresponds of course to critical points of $\mathcal{E}$. In the presence of the TT tensor $T$ however, the Euler-Lagrange shows a reaction term. Such tensors are hard to come by in general, but in dimension 4, we have encountered at the end of the previous section the Bach tensor $B$ which is a TT tensor. The following statement gives some precisions about that case. 
\begin{Th}\label{bachwill}
The equation
\be\label{bawi}
\vec{\mathcal{W}}\;=\;-\langle B,\bh\rangle\:,
\ee
is the Euler-Lagrange equation of $\mathcal{E}$ with constrained Bach energy $\int|W|^2$. 
\end{Th}
$\hfill\square$\\

\noindent
As the Bach energy is conformally invariant, solutions of (\ref{bawi}) may be called ``conformally constrained" critical points of $\mathcal{E}$. I have not been able to produce any explicit solution for this equation. In \cite{Ber2}, it is shown that limit of Palais-Smale sequences satisfy (\ref{bawo}) in the sense of distributions, but I do not know any explicit circumstance in which (\ref{bawo}) is satisfied. Naive choices, such as parallel mean curvature submanifolds, i.e. with $\pro\nabla^i\bH=\vec{0}$ for all $i$, do give rise to a TT-tensor of the right order, namely $T=|\bH|^2\bH\cdot\bh_0$, but it does not satisfy (\ref{bawo}). Another choice imposes the Chen biharmonic condition $\Delta\bH=\vec{0}$. Chen biharmonic submanifolds have been extensively considered in the literature, notably in the context of the Chen conjecture stating that $\Delta\bH=\vec{0}$ is equivalent to $\bH=\vec{0}$. \cite{Che}. An elementary computation gives that
$$
\Delta \bH\;=\;\Delta_\perp\bH-(\bH\cdot\bh^{ij})\bh_{ij}+2\nabla_i\big(|\bH|^2g^{ij}-\bH\cdot\bh^{ij}\big)\nabla_j\bp\:.
$$
For a biharmonic submanifold, the tensor $T^{ij}:=|\bH|^2g^{ij}-\bH\cdot\bh^{ij}$ is thus a TT-tensor, but it does not have the right order to serve in (\ref{bawo}). Devising a suitable nontrivial example of a solution to (\ref{bawo}) remains unclear to me. In \cite{Ber2}, it is shown that a local Palais-Smale sequence for the energy $\mathcal{E}_A$ converges to a solution of (\ref{bawo}) for some $TT$ tensor $T$. \\

Let us return to the $\bL$-system (\ref{sysL20}), and demand that $\bL$ be a closed form. It may not be immediately manifest, but as we will show, the $\bL$-system is in divergence form. For example, the first equation may be expressed as
\be\label{ojo1}
\nabla_i\big(\bL^{ij}\cdot\nabla_j\bp-\nabla^i|\bH|^2\big)\;=\;0\:.
\ee
This equation can be integrated once and yield per the Poincar\'e lemma a 2-form $S$ with
$$
d^\star S\;=\;\bL\stackrel{\!\!\!\cdot}{\res} d\bp-d|\bH|^2\:.
$$
Since $\bL$ is closed, we may impose $dS=2\bL\stackrel{\cdot}{\wedge}d\bp$. The notation used is clarified in the Appendix.\\
Similarly -- and this is admittedly non-trivial to see -- the second equation in (\ref{sysL20}) is in divergence form too and may be integrated once to yield the existence of a 2-vector-valued 2-form $\bR$ such that
$$
d^\star \bR\;=\;\bL\stackrel{\!\!\!\wedge}{\res} d\bp+\bU\:,
$$
for some 2-vector-valued 1-form $\bU$. We again set $d\bR=2\bL\stackrel{\wedge}{\wedge}d\bp$. Naturally, neither $S$ nor $\bR$ are unique. Various conditions suitable for analysis may be further imposed so as to fix uniqueness. \\
Just as in the 2-dimensional Willmore problem, what comes as a surprise is that the potentials $S$ and $\bR$ are related by a rather simple system. The following statement is the analogue of (\ref{sysUV1})\footnote{The reader is invited to consult the Appendix for the notation used.}.
\begin{Th}\label{ThSR}
When the $\bL$-system is satisfied for a 1-vector-valued 2-form $\bL\in\Lambda^2(\Lambda^1)$, then there exist $S\in\Lambda^2(\Lambda^0)$ and $\bR\in\Lambda^2(\Lambda^2)$ satisfying
\begin{equation}\label{sysSR20}
\left\{\begin{array}{rcl}
-3d^\star(\bR-\bX\stackrel{\wedge}{\wedge}d\bp)+3d\langle\bX\stackrel{\wedge}{,}d\bp\rangle&=&\vec{\eta}\stackrel{\!\!\!\bullet}{\res}d^\star\bR+d\bR\stackrel{\!\!\!\bullet}{\res}\vec{\eta}+\vec{\eta}\res d^\star S-dS\res\vec{\eta}+\mathcal{O}\big(|\bh|^3+|\bh||\pro d\bH|\big)\\
3d^\star S&=&\vec{\eta}\stackrel{\!\!\!\cdot}{\res}d^\star\bR-d\bR\stackrel{\!\!\!\cdot}{\res}\vec{\eta}+\mathcal{O}\big(|\bh|^3+|\bH||\pro d\bH|\big)\:,
\end{array}\right. 
\end{equation}
where $\vec{\eta}:=\dfrac{1}{2}d\bp\stackrel{\wedge}{\wedge}d\bp$, and $\bX$ is the 1-vector-valued 1-form with coordinate
$$
\bX^i\;=\;\nabla^i\bH-2(|\bH|^2g^{ij}-\bh^{ij}\cdot\bH)\nabla_j\bp
$$
\end{Th}
$\hfill\square$\\

\noindent
This system, henceforth referred to as the ``$(S,\bR)$-system" is of course more complicated than its 2nd order analogue, as it involves 2-forms rather than real-valued functions. In addition, the four dimensional version (\ref{sysSR20}) does not close: it bears perturbative terms. That being said, what is of importance is that the perturbative terms involve neither $\bL$ nor $\bU$, and in particular nor $\Delta_\perp\bH$, which are present at the start of the argument. As such (\ref{sysSR20}) turns out to be quite amenable to analytical considerations, under the right circumstances of course. There persists in the (\ref{sysSR20}) the term $\bX$, which does involve $d\bH$, but it appears in such a way that it will be analytically manageable. \\

Analogously to what happened in dimension two, the question of the return to geometric information arises. Just as in dimension two, it manifests itself in the form of an identity in divergence-form and which more or less corresponds to the invariance of the energy under special conformal transformation (e.g. inversion). More precisely, we have with the same vector-field $\bX$ as in Theorem~\ref{ThSR}\footnote{It may seem surprising that $\bX$ does not depend on the energy $\mathcal{E}_0$. This will be elucidated in the proof. }
$$
-2d^\star\bX\;=\;d^\star\big(\bR\stackrel{\res}{\res}d\bp-S\res d\bp\big)\:.
$$
Pursuing computations further, we will see that this implies
$$
2d^\star\bX\;=\;\pro d^\star\big((\langle\vec{\eta}\stackrel{\bullet}{,}\bR\rangle+\langle\vec{\eta},S\rangle)\res d\bp\big)+\mathcal{O}\big(|\bh||\bR|\big)\:.
$$
In particular, we will establish as a corollary:
\begin{Th}\label{Threturn}
With $\bR$ and $S$ as above, it holds
\begin{eqnarray*}
-\Delta_\perp\bH&=&\pro d^\star(\bq\res d\bp\big)+\mathcal{O}\big(|\bh||\bR|+|\bh|^3\big)\:,
\end{eqnarray*}
where $\bq$ is the 2-vector-valued scalar
$$
\bq\;:=\;\langle\vec{\eta}\stackrel{\bullet}{,}\bR\rangle+\langle\vec{\eta},S\rangle+3\langle\bX\stackrel{\wedge}{,}d\bp\rangle\:.
$$
\end{Th}
$\hfill\square$\\

\noindent
Finally, information about the 2-vector-valued scalar $\bq$ follows from a consequence of Theorem~\ref{ThSR}, namely:
\begin{Co}\label{COSR}
It holds
$$
\Delta \bq\;=\;d^{\star}\bv\qquad\text{with}\qquad \bv\;=\;\mathcal{O}\big(|\bh|^3+|\bh||\pro d\bH|+|\bh||\bR|+|\bh||S|\big)\:.
$$
\end{Co}
$\hfill\square$\\

The four-dimensional situation bears striking resemblance to the two-dimensional Willmore problem. Some differences arise, mostly due to the fact that the $(S,\bR)$-system does not close. Nevertheless, the various highly structured equations satisfied by critical points make it possible to obtain analytical information, under some smallness hypothesis on $\bh$. This is the main subject of \cite{Ber2}.

\setcounter{equation}{0} 
\reset

\section{Proofs of the Theorems}

We recall once and for all the well-known expressions \cite{KN} of the Riemann, Ricci, and scalar curvature tensors for a submanifold: 
$$
R^{ijkl}\;=\;\bh^{ik}\cdot\bh^{jl}-\bh^{il}\cdot\bh^{jk}\quad,\quad R^{ik}\;=\;4\bH\cdot\bh^{ik}-\bh^{ij}\cdot\bh^k_j\quad,\quad R\;=\;16|\bH|^2-|\bh|^2\:.
$$
Moreover, we recall the classical expression 
$$
|W|^2\;=\;|\text{Rm}|^2-2|\text{Ric}|^2+\dfrac{1}{3}R^2\:.
$$

\subsection{Conformally invariant energies: proof of Proposition~\ref{EC}}\label{UC}

We denote by $(h^\al)^{ij}$ the normal component of $\bh^{ij}$, with $\al\in\{1,\ldots,\text{codim}\}$. Moreover, for convenience, we let
$$
(h^2)^{ij}\;:=\;\bh^{ik}\cdot\bh_k^j\:.
$$
Also, we let 
$$
\langle \bh^4\rangle\;:=\;(\bh^{ik}\cdot\bh^{jl})(\bh_{ik}\cdot\bh_{jl})\:.
$$
A classical corollary of Simons' identity (see e.g. \cite{Smo}) asserts that
\begin{eqnarray}\label{simo0}
8\langle\bh\stackrel{\cdot}{,}\nabla^2\bH\rangle&=&\Delta |\bh|^2-2|\pro d\bh|^2+|\text{Rm}|^2+2|\text{Ric}|^2-32|\bH\cdot\bh|^2\nonumber\\
&&\hspace{1cm}+\:\big|(h^\al)^{ik}(h^\beta)^{j}_{k}-(h^\beta)^{il}(h^\al)^{j}_{l}\big|^2\:.
\end{eqnarray}
The last term is expanded as
\begin{eqnarray*}
\big|(h^\al)^{ik}(h^\beta)^{j}_{k}-(h^\beta)^{il}(h^\al)^{j}_{l}\big|^2&=&2|h^2|^2-2(\bh_{ik}\cdot\bh_{jl})(\bh^{jk}\cdot\bh^{il})\\
&=&2|h^2|^2-2(\bh_{ik}\cdot\bh_{jl})R^{jikl}-2\langle \bh^4\rangle\\
&=&2|h^2|^2+|\text{Rm}|^2-2\langle \bh^4\rangle\:.
\end{eqnarray*}
Brought into (\ref{simo0}), the latter yields
\be\label{simo00}
8\langle\bh\stackrel{\cdot}{,}\nabla^2\bH\rangle\;=\;\Delta |\bh|^2-2|\pro d\bh|^2+2|\text{Rm}|^2+2|\text{Ric}|^2-32|\bH\cdot\bh|^2+2|h^2|^2-2\langle \bh^4\rangle\:.
\ee
Note that the left-hand side may be recast as
\begin{eqnarray*}
8\langle\bh\stackrel{\cdot}{,}\nabla^2\bH\rangle\;=\;8\nabla_{ij}\big(\bH\cdot\bh^{ij}-4|\bH|^2g^{ij}\big)+16\Delta |\bH|^2-32|\pro d\bH|^2-8\langle \bH\cdot h^2\stackrel{}{,}\bh\rangle\:.
\end{eqnarray*}
Hence now, recalling that $R=16|\bH|^2-|\bh|^2$, 
\begin{eqnarray}\label{simo1}
8\nabla_{ij}\big(\bH\cdot\bh^{ij}-4|\bH|^2g^{ij}\big)+\Delta R&=&32|\pro d\bH|^2-2|\pro d\bh|^2+2|\text{Rm}|^2+2|\text{Ric}|^2-32|\bH\cdot\bh|^2\nonumber\\
&&+\:2|h^2|^2+8\langle \bH\cdot\bh\stackrel{}{,}h^2\rangle-2\langle \bh^4\rangle\:.
\end{eqnarray}
In order to simplify this expression, we first note that
\begin{eqnarray*}
(h^2)^{ij}&=&\bh^{ik}\cdot\bh_{k}^j\;\;=\;\;(\bh_0^{ik}+g^{ik}\bH)\cdot((\bh_0)_k^j+g_k^j\bH)\\
&=&(h_0^2)^{ij}+2\bH\cdot\bh^{ij}-|\bH|^2g^{ij}\:.
\end{eqnarray*}
Hence
\begin{eqnarray*}
|h^2|^2&=&|h_0^2|^2+4|\bH\cdot\bh|^2-12|\bH|^4-2|\bH|^2|\bh_0|^2+4\langle \bH\cdot\bh\stackrel{}{,}\bh_0^2\rangle\\[1ex]
&=&|h_0^2|^2+4|\bH\cdot\bh|^2-12|\bH|^4-2|\bH|^2(|\bh|^2-4|\bH|^2)\\
&&\hspace{1cm}+\:4\langle \bH\cdot\bh\stackrel{}{,}h^2-2\bH\cdot\bh+|\bH|^2g\rangle\\[1ex]
&=&|h_0^2|^2-4|\bH\cdot\bh|^2-2|\bH|^2|\bh|^2+4\langle \bH\cdot\bh\stackrel{}{,}h^2\rangle+12|\bH|^4\:,
\end{eqnarray*}
so that
\begin{eqnarray*}
4|\text{Ric}|^2&=&4|h^2|^2-32\langle \bH\cdot\bh\stackrel{}{,}h^2\rangle+64|\bH\cdot\bh|^2\\[.5ex]
&=&6|h^2|^2-2|h^2|^2-32\langle \bH\cdot\bh\stackrel{}{,}h^2\rangle+64|\bH\cdot\bh|^2\\[.5ex]
&=&-\,2|h^2|^2-8\langle \bH\cdot\bh\stackrel{}{,}h^2\rangle+6|h_0^2|^2+40|\bH\cdot\bh|^2-12|\bH|^2|\bh|^2+72|\bH|^4\:.
\end{eqnarray*}
Introduced into (\ref{simo1}), the latter shows that
\begin{eqnarray}\label{simob1}
8\nabla_{ij}\big(\bH\cdot\bh^{ij}-4|\bH|^2g^{ij}\big)+\Delta R&=&32|\pro d\bH|^2-2|\pro d\bh|^2+2|\text{Rm}|^2-2|\text{Ric}|^2+6|h_0^2|^2\nonumber\\
&&\hspace{.5cm}+\,8|\bH\cdot\bh|^2-12|\bH|^2|\bh|^2+72|\bH|^4-2\langle \bh^4\rangle\:.
\end{eqnarray}

%
%
%
To further simplify, we compute
\begin{eqnarray*}
\bh^{ik}\cdot\bh^{jl}&=&\bh_0^{ik}\cdot\bh_0^{jl}+\bH\cdot(\bh^{jl}g^{ik}+\bh^{ik}g^{jl})-|\bH|^2g^{ik}g^{jl}\:,
\end{eqnarray*}
from which it easily follows that
\begin{eqnarray*}
\langle \bh^4\rangle&=&\langle \bh_0^4\rangle+8|\bH\cdot\bh|^2-16|\bH|^4\:.
\end{eqnarray*}
Introducing this into (\ref{simob1}) gives
\begin{eqnarray}\label{simo3}
8\nabla_{ij}\big(\bH\cdot\bh^{ij}-4|\bH|^2g^{ij}\big)+\Delta R&=&32|\pro d\bH|^2-2|\pro d\bh|^2+2|\text{Rm}|^2-2|\text{Ric}|^2\nonumber\\
&&\hspace{1cm}-\:8|\bH\cdot\bh|^2-12|\bH|^2|\bh|^2+104|\bH|^4\nonumber\\
&&\hspace{2cm}+\:6|h_0^2|^2-2\langle \bh_0^4\rangle\:.
\end{eqnarray}
We note next that
\begin{eqnarray*}
2|\text{Rm}|^2-2|\text{Ric}|^2&=&3\Big(|\text{Rm}|^2-2|\text{Ric}|^2+\dfrac{1}{3}R^2\Big)-\big(|\text{Rm}|^2-4|\text{Ric}|^2+R^2\big)\\
&=&3|\text{W}|^2-c(\Sigma)\:,
\end{eqnarray*}
where $c(\Sigma)$ denotes the integrand in the Chern-Gauss-Bonnet topological invariant. The latter is now injected into (\ref{simo3}) to yield
\begin{eqnarray}\label{simo4}
&&8\nabla_{ij}\big(\bH\cdot\bh^{ij}-4|\bH|^2g^{ij}\big)+\Delta R+c(\Sigma)\nonumber\\
&&\hspace{3cm}=\:32\big(|\pro d\bH|^2-|\bH\cdot\bh|^2+7|\bH|^4)\nonumber\\
&&\hspace{4cm}-\:2\big(|\pro d\bh|^2-12|\bH\cdot\bh|^2+6|\bH|^2|\bh|^2+60|\bH|^4\big)\nonumber\\
&&\hspace{5cm}+\:6|h_0^2|^2-2\langle \bh_0^4\rangle+3|W|^2\:.
\end{eqnarray}
Therefore, upon integration over a closed manifold $\Sigma$ without boundary, we arrive at
$$
16\mathcal{E}_A-\mathcal{E}_C+\mathcal{E}_0\;=\;16\pi^2\chi(\Sigma)-\dfrac{3}{2}\int_\Sigma|W|^2d\text{vol}_g\:,
$$
where
$$
\mathcal{E}_0\;:=\;\int_\Sigma\big(3|h_0^2|^2-\langle \bh_0^4\rangle\big)d\text{vol}_g\:,
$$
and
$$
\mathcal{E}_C\;:=\;\int_\Sigma\big(|\pro d\bh|^2-12|\bH\cdot\bh|^2+6|\bH|^2|\bh|^2+60|\bH|^4\big)d\text{vol}_g\:.
$$
As $\mathcal{E}_A$, $\int|W|^2$, and $\mathcal{E}_0$ are all conformally invariant, it follows that $\mathcal{E}_C$ is likewise conformally invariant. This is the claim of Proposition~\ref{EC}.

\subsection{The first conservation law: proof of Theorem~\ref{cons}}

Consider a conformally invariant energy of the form
$$
\mathcal{E}\;=\;\mathcal{E}_A+\mathcal{E}_0\:,
$$
where $\mathcal{E}_A$ is as in the previous section, while $\mathcal{E}_0$ is any energy with integrand of order $\mathcal{O}(|\bh_0|^4)$, see Subsection~\ref{derder}. As the Euler-Lagrange equation is quite unmanageable as written in (\ref{horrib}), the first thing to do is to express it in divergence-form. This can be accomplished by evaluating the Noether field corresponding to translation invariance. Per the considerations in Section~\ref{Noether}, a critical point of $\mathcal{E}$ satisfies a divergence-form equation
$$
d^\star\bV\;=\;\vec{0}\:,
$$
where the 1-vector-valued 1-form $\bV$ satisfies
\be\label{defV}
\bV^j\;=\;G^{ij}\nabla_i\bp-\pro\nabla_i\bF^{ij}+\bC^j\:.
\ee
Exact expressions for $G$, $\bF$, and $\bC$ can be found in the Appendix. For the moment, we record that
\be\label{GAFO}
\left\{\begin{array}{lcl}
G^{ij}&=&\dfrac{1}{2}\bh^{ij}\cdot\Delta_\perp\bH+ |\pro\nabla\bH|^2g^{ij}-2\nabla^j\bH\cdot\pro\nabla^i\bH+\mathcal{O}(|\bh|^4)\\[1ex]
\bF^{ij}&=&-\dfrac{1}{2}\Delta_\perp\bH\,g^{ij}+\mathcal{O}(|\bh|^3) \\[1ex]
\bC^j&=&-2(\bh^{ij}\cdot\nabla_i\bH)\bH+2(\bH\cdot\bh^{ij})\pro\nabla_i\bH\:.\end{array}\right.
\ee
Since
$$
-\pro\nabla^j\Delta_\perp\bH\;=\;-\nabla^j\Delta_\perp\bH+\pi_T\nabla^j\Delta_\perp\bH\;=\;-\nabla^j\Delta_\perp\bH-\big(\bh^{ij}\cdot\Delta_\perp\bH\big)\nabla_i\bp\:,
$$
we can recast (\ref{defV}) in the form
\be\label{defVV}
\bV^j\;=\;\dfrac{1}{2}\nabla^j\Delta_\perp\bH+G_1^{ij}\nabla_i\bp-\pro\nabla_i\bF_1^{ij}+\mathcal{O}\big(|\bh|^2|\pro d\bH|\big)\:,
\ee
with
\be\label{GAFB}
\left\{\begin{array}{lcl}
G_1^{ij}&=&\bh^{ij}\cdot\Delta_\perp\bH+ |\pro\nabla\bH|^2g^{ij}-2\nabla^j\bH\cdot\pro\nabla^i\bH+\mathcal{O}(|\bh|^4)\\[1ex]
\bF_1^{ij}&=&\mathcal{O}(|\bh|^3)\:.\end{array}\right.
\ee
A simple application of the Codazzi identity gives
\begin{eqnarray*}
&&\nabla_k(\bh^{ij}\cdot\nabla^k\bH)-\nabla^i(\bh^{kj}\cdot\nabla_k\bH)+\nabla_k(\bh^{kj}\cdot\nabla^i\bH)\\[1ex]
&=&\bh^{ij}\cdot\Delta_\perp\bH+\pro(\nabla^k\bh^{ij}-\nabla^i\bh^{kj})\cdot\nabla_k\bH+4\nabla^j\bH\cdot\pro\nabla^i\bH+\bh^{j}_k\cdot\big(\nabla^k\pro\nabla^i\bH-\nabla^i\pro\nabla^k\bH\big)\\[1ex]
&=&\bh^{ij}\cdot\Delta_\perp\bH+\mathcal{O}\big(|\pro d\bH|^2\big)-\bh^{j}_k\cdot\big(\nabla^k\pi_T\nabla^i\bH-\nabla^i\pi_T\nabla^k\bH\big)\\[1ex]
&=&\bh^{ij}\cdot\Delta_\perp\bH+\mathcal{O}\big(|\pro d\bH|^2\big)+\bh^{j}_k\cdot\Big(\nabla^k\big((\bH\cdot\bh^{ij})\nabla_j\bp\big)-\nabla^i\big((\bH\cdot\bh^{kj})\nabla_j\bp\big)\Big)\\[1ex]
&=&\bh^{ij}\cdot\Delta_\perp\bH+\mathcal{O}\big(|\pro d\bH|^2+|\bh|^4\big)\:.
\end{eqnarray*}
Accordingly, 
\begin{eqnarray*}
(\bh^{ij}\cdot\Delta_\perp\bH)\nabla_i\bp&=&\nabla_k\big((\bh^{ij}\cdot\nabla^k\bH+\bh^{kj}\cdot\nabla^i\bH)\nabla_i\bp-(\bh^{ij}\cdot\nabla_i\bH)\nabla^k\bp\big)+\mathcal{O}\big(|\pro d\bH|^2+|\bh|^4\big)\:.
\end{eqnarray*}
Introducing this into (\ref{defVV}) and using the fact that $\pro\nabla_i\bF_1^{ij}=\nabla_i\bF_1^{ij}+\mathcal{O}(|\bF_1||\bh|)$, since $\bF_1$ is normal valued, we arrive at
\bes
\bV^j\;=\;\dfrac{1}{2}\nabla^j\Delta_\perp\bH+\nabla_kA^{kj}+\mathcal{O}\big(|\bh|^2|\pro d\bH|+|\pro d\bH|^2+|\bh|^4\big)\:,
\ees
with
$$
A^{kj}\;:=\;(\bh^{ij}\cdot\nabla^k\bH+\bh^{kj}\cdot\nabla^i\bH)\nabla_i\bp-(\bh^{ij}\cdot\nabla_i\bH)\nabla^k\bp-\bF_1^{kj}\;=\;\mathcal{O}\big(|\bh||\pro d\bH|+|\bh|^3\big)\:.
$$
This completes the proof of Theorem~\ref{cons}.

\subsection{The $\bL$-system: proofs of Theorems~\ref{ThsysL},~\ref{ThPS}, and~\ref{bachwill}}

Let us return to (\ref{GAFO}). The tensors $G$, $\bF$, and $\bC$ are studied at length in Section~\ref{Noether}, and we summarise here their properties of immediate interest. 

\begin{Prop}\label{prop1}
$G$ and $\bF$ are symmetric two-tensors, while $\bF$ and $\bC$ are normal-valued. In addition, the following identities hold:
\begin{itemize}
\item[(i)] $G^i_i\;=\;\Delta|\bH|^2$
\item[(ii)] $\bh^m_j\cdot\big(-\nabla_i\bF^{ij}+\bC^j\big)\;=\;\nabla_iG^{im}\qquad\forall\:m$
\item[(iii)] $(-\pro\nabla_i\bF^{ij}+\bC^j)\wedge\nabla_j\bp\;=\;\nabla_i\big(-\bF^{ij}\wedge\nabla_j\bp+2\bH\wedge\pro\nabla^i\bH\big)\:.$
\end{itemize}
\end{Prop}
{\bf Proof.} The first item is Remark~\ref{rem0}. The second item is obtained combining Lemmata~\ref{lem0} and ~\ref{lem3}, as well as Corollary~\ref{cor1}. The third item follows directly from Lemma~\ref{lem4} and Corollary~\ref{cor2}.  

$\hfill\blacksquare$\\

From the Poincar\'e Lemma, we know there exists a 1-vector-valued 2-form $\bL_0$ satisfying
\be\label{eqq1}
(d^\star\bL_0)^j\;=\;\bV^j\;=\;G^{ij}\nabla_i\bp-\pro\nabla_i\bF^{ij}+\bC^j\:.
\ee
Using (i) from Proposition~\ref{prop1} shows that
\be\label{L1}
(d^\star\bL_0)^j\cdot\nabla_j\bp\;=\;G^i_i\;=\;\Delta|\bH|^2\:.
\ee
Moreover, since $G$ is symmetric, we have
\be\label{L2}
(d^\star\bL_0)^{j}\wedge\nabla_j\bp\;=\;(-\pro\nabla_i\bF^{ij}+\bC^j)\wedge\nabla_j\bp\:.
\ee
Together, (\ref{L1})-(\ref{L2}) form the so-called {\it $\bL$-system}, which is satisfied for a critical point of the energy $\mathcal{E}$. The following natural question arises: if there exists $\bL\in\Lambda^2(\Lambda^1)$ such that the $\bL$-system is satisfied, does it follow that we are dealing with a critical point of $\mathcal{E}$? The answer is ``not quite", as we will show next. \\

We suppose that some form $\bL_0\in\Lambda^2(\Lambda^1)$ satisfies the system
\begin{equation}\label{sysL}
\left\{\begin{array}{lcl}(d^\star\bL_0)^j\cdot\nabla_j\bp&=&d^\star d|\bH|^2\\[1ex]
(d^\star\bL_0)^{j}\wedge\nabla_j\bp&=&-\pro\nabla_i\bF^{ij}\wedge\nabla_j\bp+\bC^j\wedge\nabla_j\bp\:.
\end{array}\right.
\end{equation}
Let us first define the 2-form $\bL:=\bL_0-\dfrac{1}{3}d(|\bH|^2d\bp)$, so that
\begin{eqnarray}\label{eqq3}
(d^\star\bL)^j\cdot\nabla_j\bp&=&(d^\star\bL_0)^j\cdot\nabla_j\bp-\dfrac{1}{3}\nabla_i\big(\nabla^i|\bH|^2\nabla^j\bp-\nabla^j|\bH|^2\nabla^i\bp)\cdot\nabla_j\bp\nonumber\\
&=&(d^\star\bL_0)^j\cdot\nabla_j\bp-\dfrac{1}{3}\big(\Delta|\bH|^2\nabla^j\bp-\nabla^{ij}|\bH|^2\nabla^i\bp)\cdot\nabla_j\bp\nonumber\\
&\stackrel{\text{(\ref{L1})}}{=}&\Delta|\bH|^2-\Delta|\bH|^2\;\;=\;\;0\:.
\end{eqnarray}
Moreover, using the Codazzi identity:
\begin{eqnarray}\label{eqq4}
(d^\star\bL)^j\wedge\nabla_j\bp&=&(d^\star\bL_0)^j\wedge\nabla_j\bp-\dfrac{1}{3}\nabla_i\big(\nabla^i|\bH|^2\nabla^j\bp-\nabla^j|\bH|^2\nabla^i\bp)\wedge\nabla_j\bp\nonumber\\
&=&(d^\star\bL_0)^j\wedge\nabla_j\bp-\dfrac{1}{3}\big(\bh^j_i\nabla^i|\bH|^2-4\bH\nabla^j|\bH|^2)\wedge\nabla_j\bp\nonumber\\
&\stackrel{\text{(\ref{L2})}}{=}&(-\pro\nabla_i\bF^{ij}+\bC^j)\wedge\nabla_j\bp-\dfrac{1}{3}\pro\nabla_i\vec{f}^{\,ij}\wedge\nabla_j\bp\:,
\end{eqnarray}
where
\be\label{deff}
\vec{f}^{\,ij}\;:=\;|\bH|^2\big(\bh^{ij}-4\bH g^{ij}\big)\:.
\ee
Equations (\ref{eqq3})-(\ref{eqq4}) help us recast (\ref{sysL}) in the form
\begin{equation}\label{sysL2}
\left\{\begin{array}{lcl}(d^\star\bL)^j\cdot\nabla_j\bp&=&0\\[1ex]
(d^\star\bL)^{j}\wedge\nabla_j\bp&=&-\pro\nabla_i\vec{\mathcal{F}}^{ij}\wedge\nabla_j\bp+\bC^j\wedge\nabla_j\bp\:,
\end{array}\right.
\end{equation}
where for notational convenience we have set
$$
\vec{\mathcal{F}}\;:=\;\bF+\dfrac{1}{3}\vec{f}\:.
$$
It is to be noted that $\vec{\mathcal{F}}$ satisfies (iii) from Proposition~\ref{prop1}. That is rather plain to see, since clearly $\vec{f}^{\,ij}\wedge\bh_{ij}=\vec{0}$. However, $\vec{\mathcal{F}}$ fails to satisfy (ii). Indeed, since $\pro\nabla_i(\bh^{ij}-4\bH g^{ij})=\vec{0}$ by Codazzi, it comes
\begin{eqnarray*}
\bh^m_j\cdot\nabla_i\vec{f}^{\,ij}&=&\bh^m_j\cdot(\bh^{ij}-4\bH g^{ij})\nabla_i|\bH|^2\;=\;-R^{mi}\nabla_i|\bH|^2\:,
\end{eqnarray*}
where the reader is invited to recall the extrinsic expression of the Ricci tensor $R^{im}$. From this observation and Proposition~\ref{prop1}, we see that
\be\label{eqq6}
\bh^m_j\cdot\big(-\nabla_i\vec{\mathcal{F}}^{ij}+\bC^j\big)\;=\;\nabla_iG^{im}+\dfrac{1}{3}R^{im}\nabla_i|\bH|^2\:.
\ee
Define next the symmetric two-tensor
$$
a^{im}\;:=\;\nabla^{im}|\bH|^2-g^{im}\Delta|\bH|^2\:.
$$
We note that
$$
\nabla_ia^{im}\;=\;(\nabla^{ij}-\nabla^{ji})\nabla_i|\bH|^2\;=\;R^{im}\nabla_i|\bH|^2\:.
$$
Combining this to (\ref{eqq6}), it comes
\be\label{eqq7}
\bh^m_j\cdot\big(-\nabla_i\vec{\mathcal{F}}^{ij}+\bC^j\big)\;=\;\nabla_j\Big(G^{jm}+\dfrac{1}{3}a^{jm}\Big)\:.
\ee

We next introduce the projection decomposition $(d^\star\bL)^{j}=:A^{ji}\nabla_i\bp+\bB^j$, where $\bB$ is normal-valued. Introducing this into (\ref{sysL2}), we find
\begin{equation*}
\left\{\begin{array}{rcl}A^i_i&=&0\\[1ex]
\dfrac{1}{2}(A^{ji}-A^{ij})\vec{\eta}_{ij}+\bB^j\wedge\nabla_j\bp&=&-\pro\nabla_i\vec{\mathcal{F}}^{ij}\wedge\nabla_j\bp+\bC^j\wedge\nabla_j\bp\:.
\end{array}\right.
\end{equation*}
In other words, $A$ is symmetric and traceless, while the normal part satisfies
$$
\bB^j\;=\;-\pro\nabla_i\vec{\mathcal{F}}^{ij}+\bC^j\:.
$$
Note that the tangential projection satisfies
$$
\vec{0}\;=\;\pi_Td^\star d^\star\bL\;=\;\pi_T\nabla_j(A^{ji}\nabla_i\bp+\vec{B}^j)\;=\;\big(\nabla_jA^{ij}-\vec{B}^{j}\cdot\bh_j^i\big)\nabla_i\bp\:,
$$
where we have used that $\vec{B}$ is normal-valued. It follows that
$$
\nabla_jA^{jm}\;=\;\bB^j\cdot\bh_j^m\;=\;\bh^m_j\cdot\big(-\nabla_i\vec{\mathcal{F}}^{ij}+\bC^j\big)\;\stackrel{\text{(\ref{eqq7})}}{=}\;\nabla_j\Big(G^{jm}+\dfrac{1}{3}a^{jm}\Big)\:.
$$
Accordingly, the symmetric 2-tensor $T:=A-G-\frac{1}{3}a$ is divergence-free. Moreover, from (i) in Proposition~\ref{prop1} and the tracelessness of $A$, we easily verify that $T$ is also traceless, thereby making $T$ into a traceless transverse tensor (abbreviated {\it TT-tensor}). \\[1.5ex]
Summarising our findings so far, it holds
\begin{eqnarray}\label{eqq8}
(d^\star\bL)^j&=&\Big(T^{ij}+G^{ij}+\dfrac{1}{3}a^{ij}\Big)\nabla_i\bp-\pro\nabla_k\vec{\mathcal{F}}^{kj}+\bC^j\nonumber\\
&\stackrel{\text{(\ref{defV})}}{=}&\bV^j+\Big(T^{ij}+\dfrac{1}{3}a^{ij}\Big)\nabla_i\bp-\dfrac{1}{3}\pro\nabla_i\vec{f}^{\,ij}\:.
\end{eqnarray}
To further proceed, we need the following
\begin{Lm}
We have
$$
\nabla_i\big(a^{ij}\nabla_j\bp-\pro\nabla_j\vec{f}^{\,ij}\big)\;=\;\vec{0}\:.
$$
\end{Lm}
{\bf Proof.} We have already observed that 
\be\label{eqq10}
\pro\nabla_j\vec{f}^{\,ij}\;=\;(\bh^{ij}-4\bH g^{ij})\nabla_j|\bH|^2\:,
\ee
so that, recycling some of our previous derivations, we obtain
\begin{eqnarray*}
\nabla_i\pro\nabla_j\vec{f}^{\,ij}&=&(\bh^{ij}-4\bH g^{ij})\nabla_{ij}|\bH|^2-(\bh^{ij}-4\bH g^{ij})\cdot\bh_i^k\nabla_j|\bH|^2\nabla_k\bp\\[1ex]
&=&(\bh^{ij}-4\bH g^{ij})\nabla_{ij}|\bH|^2+R^{ji}\nabla_j|\bH|^2\nabla_i\bp\\
&=&a_{ij}\bh^{ij}+\nabla_ja^{ij}\nabla_i\bp\\
&=&\nabla_j(a^{ij}\nabla_i\bp)\:.
\end{eqnarray*}
This confirms that the vector-field $a^{ij}\nabla_j\bp-\pro\nabla_j\vec{f}^{\,ij}$ is co-closed. \\

$\hfill\blacksquare$\\

Using this lemma and the fact that $T$ is divergence-free we may now apply $d^\star$ to both sides of (\ref{eqq8}) so as to arrive at
$$
\vec{0}\;=\;\nabla_j\bV^j+\langle T,\bh\rangle\:.
$$
In other words, whereas critical points of $\mathcal{E}$ satisfy $d^\star\bV=\vec{0}$ and the $\bL$-system, if only the $\bL$-system is satisfied, then the perturbed equation $d^\star\bV=-\langle T,\bh\rangle$ holds for some TT-tensor $T$. Conversely, if the equation $d^\star\bV=-\langle T,\bh\rangle$ holds for some TT-tensor $T$, noting that 
$$
\langle T,\bh\rangle\;=\;d^\star(T^{ij}\nabla_i\bp)\:,
$$
we deduce the existence of a two-form $\bL$ such that
$$
(d^\star\bL)^j\;=\;\bV^j+T^{ij}\nabla_i\bp\:.
$$
Since $T$ is a symmetric $TT$-tensor, it immediately follows that
$$
(d^\star\bL)^j\cdot\nabla_j\bp\;=\;\bV^j\cdot\nabla_j\bp\;=\;G^i_i\:,
$$
as well as
$$
(d^\star\bL)^j\wedge\nabla_j\bp\;=\;\bV^j\wedge\nabla_j\bp\;=\;-\pro\nabla_i\bF^{ij}\wedge\nabla_j\bp+\bC^j\wedge\nabla_j\bp\:.
$$
Therefore, we recover the $\bL$-system in that case. We are of course free to demand that $\bL$ be closed. This completes the proof of Theorem~\ref{ThPS}.\\

As an example of the $\bL$-system, we consider the case when $T$ is chosen to be the Bach tensor $B$. We claim that the Euler-Lagrange equation corresponds to varying $\mathcal{E}$ while constraining the Bach energy $\int|W|^2$. To verify this, we need to compute the Noether field associated with the Bach energy (understood extrinsically). This is done in the Appendix, and one verifies that
\be\label{ojo2}
\bV^b_{|W|^2}\;=\;-4P_{ij}W^{biaj}\nabla_a\bp-8\pro\nabla_i\big(\bh_{aj}W^{abji}\big)\:.
\ee
For the next result, we need to recall the Cotton form 
$$
\vec{\ell}^{ab}\;:=\;C^{abc}\nabla_c\bp\;\equiv\;(\nabla^aP^{bc}-\nabla^bP^{ac})\nabla_c\bp\;\equiv\;2\nabla_iW^{icab}\nabla_c\bp\:.
$$
Using the Cotton-Weyl equation as well as the symmetries of the Weyl tensor and the Codazzi equation gives easily
\begin{eqnarray*}
(d^\star\vec{\ell})^b&=&\nabla_aC^{abc}\nabla_c\bp+C^{abc}\bh_{ca}\\
&=&\big(B^{bc}-P_{ij}W^{bicj}\big)\nabla_c\bp+2\bh_{ca}\nabla_iW^{icab}\\
&=&\big(B^{ba}-P_{ij}W^{biaj}\big)\nabla_a\bp-2\pro \nabla_i\big(\bh_{ja}W^{abji}\big)\:.
\end{eqnarray*}
Combined with (\ref{ojo2}), the latter yields that
\be
-\dfrac{1}{4}\bV^i_{|W|^2}+(d^\star\vec{\ell})^i\;=\;B^{ij}\nabla_j\bp\:.
\ee
As the Bach tensor is divergence-free, it follows that extrinsic critical points of the Bach energy satisfy
$$
\dfrac{1}{4}d^\star\bV_{|W|^2}\;=\;-\langle B,\bh\rangle\:.
$$
In particular, when the variation of $\mathcal{E}$ is done with the Bach energy held fixed, then one reaches the perturbed equation
$$
\vec{\mathcal{W}}\;=\;-\langle B,\bh\rangle\:.
$$
This is the content of Theorem~\ref{bachwill}.

\subsection{The $(S,\bR)$-system: proof of Theorem~\ref{ThSR} and Corollary~\ref{COSR}}

Let us return to (\ref{eqq1}). First of all, let us declare, since we have the freedom to do so, that $\bL_0$ is a closed form. As $\bC$ is normal-valued, we can use the first item in Proposition~\ref{prop1} to obtain\footnote{note that $\bL_0^{ij}\cdot\nabla_{ij}\bp=0$, since $\bL_0$ is anti-symmetric.}
\be\label{eqq9}
\nabla_i\big(\bL_0^{ij}\cdot\nabla_j\bp\big)\;=\;\nabla_i\bL_0^{ij}\cdot\nabla_j\bp\;=\;G^i_i\;=\;\nabla_i\nabla^i|\bH|^2\:. 
\ee
Setting $\bL:=\bL_0-\dfrac{1}{3}d|\bH|^2\wedge d\bp$, it comes
$$
\bL^{ij}\cdot\nabla_j\bp\;=\;\bL^{ij}_0\cdot\nabla_j\bp-\nabla^i|\bH|^2\:.
$$
Note that $\bL$ is closed. Now (\ref{eqq9}) yields
$$
\nabla_i\big(\bL^{ij}\cdot\nabla_j\bp\big)\;=\;0\:.
$$
From the Poincar\'e Lemma, we deduce the existence of a scalar-valued two-form $S$ satisfying\footnote{The reader is invited to consult Section~\ref{nota} for the notational conventions applying in this paper.}
\be\label{defS}
d^\star S\;=\;\bL\stackrel{\!\!\!\cdot}{\res}d\bp\:.
\ee
For reasons that will become clear in the sequel, we ask that
\be\label{deffS}
dS\;=\;2\bL\stackrel{\cdot}{\wedge}d\bp\:,
\ee
which is indeed a closed form, since $\bL$ is. \\

On the other hand, the second equation in (\ref{sysL2}) gives
\begin{eqnarray*}
\nabla_i\big(\bL^{ij}\wedge\nabla_j\bp\big)\;=\;-\pro\nabla_i\vec{\mathcal{F}}^{ij}\wedge\nabla_j\bp+\bC^j\wedge\nabla_j\bp\:.
\end{eqnarray*}
Recall that $\vec{\mathcal{F}}=\bF+\frac{1}{3}\vec{f}^{\,ij}$, with $\vec{f}$ given by (\ref{deff}). According to the latter and the last item in Proposition~\ref{prop1}, we find 
\begin{eqnarray}\label{eqq11}
\nabla_i\big(\bL^{ij}\wedge\nabla_j\bp\big)\;=\;\nabla_i\bU^i\:,
\end{eqnarray}
where we have set
\be\label{defU}
\bU^i\;:=\;-\bF^{ij}\wedge\nabla_j\bp+2\bH\wedge\pro\nabla^i\bH-\dfrac{1}{3}\vec{f}^{\,ij}\wedge\nabla_i\bp\:.
\ee
\begin{Rm}
For our future use, we record the following fact:
\begin{eqnarray}\label{huma3}
\bU\stackrel{\res}{\res}d\bp&=&\bF^i_i+\dfrac{1}{3}\vec{f}^i_i\;\;=\;\;-2\Delta_\perp\bH-2\langle\bH\cdot\bh,\bh\rangle+16|\bH|^2\bH\:.
\end{eqnarray}
We have used (\ref{huma1})-(\ref{huma2}) from the Appendix. This expression may be further simplified owing to an observation originally made by Guven \cite{Guv}. Upon varying the energy $\int|\bH|^2$, one finds a generic identity which holds in dimension 4\footnote{see also Lemma~\ref{lem4}.}:
\be\label{yab}
\Delta_\perp\bH+\langle\bH\cdot\bh,\bh\rangle-8|\bH|^2\bH\;=\;d^\star\bX\:,
\ee
where $\bX$ is the 1-vector-valued 1-form 
$$
\bX^i\;:=\;\nabla^i\bH-2(|\bH|^2g^{ij}-\bH\cdot\bh^{ij})\nabla_j\bp\:.
$$
We may now recast (\ref{huma3}) in the form
\be\label{huma4}
\bU\stackrel{\res}{\res}d\bp\;=\;-2d^\star\bX\:.
\ee
\end{Rm}
$\hfill\square$

We have used the already encountered observation:
$$
\pro\nabla_i\vec{f}^{\,ij}\wedge\nabla_j\bp\;=\;\nabla_i\big(\vec{f}^{\,ij}\wedge\nabla_j\bp)\:.
$$
We deduce from (\ref{eqq11}) the existence of a two-vector-valued two-form $\bR$ satisfying
\be\label{defR}
d^\star\bR\;=\;\bL\stackrel{\!\!\!\wedge}{\res}d\bp+\bU\:.
\ee
For future purposes, we set
\be\label{deffR}
d\bR\;=\;2\bL\stackrel{\wedge}{\wedge}d\bp\:.
\ee
Again, this is permitted since $\bL$ is a closed form. \\

Our aim in this section is to uncover structural identities linking $S$ and $\vec{R}$ together. In order to do so, we employ some generic results. Throughout the rest of this section, we let $\vec{\eta}:=\frac{1}{2}d\bp\stackrel{\wedge}{\wedge}d\bp$. 

\begin{Prop}\label{propito}
Let $\vec{\ell}\in\Lambda^2(\Lambda^1)$, and define accordingly
$$
A\;:=\;\vec{\ell\,}\stackrel{\!\!\!\cdot}{\res}d\bp\quad,\quad 
B\;:=\;2\vec{\ell\,}\stackrel{\cdot}{\wedge}d\bp\quad,\quad
\bCC\;:=\;\vec{\ell\,}\stackrel{\!\!\!\wedge}{\res}d\bp\quad,\quad
\bD\;:=\;2\vec{\ell\,}\stackrel{\wedge}{\wedge}d\bp\:.
$$
The following identities hold
\be\label{clarisse}
 -\,3\bCC\;=\;\vec{\eta}\stackrel{\!\!\!\bullet}{\res}\bCC+\bD\stackrel{\!\!\!\bullet}{\res}\vec{\eta}+\vec{\eta}\res A-B\res\vec{\eta}\:.
 \ee
 and
 \be\label{clarisse2}
3A\;=\;\vec{\eta}\stackrel{\!\!\!\cdot}{\res}\bCC-\bD\stackrel{\!\!\!\cdot}{\res}\vec{\eta}\:.
 \ee
\end{Prop}
{\bf Proof.}
We will first show (\ref{clarisse}). To this end, we compute:
\begin{eqnarray}\label{ahah1}
&&(\vec{\ell}^{\,ij}\wedge\nabla^k\bp)\bullet(\nabla_a\bp\wedge\nabla_b\bp)\nonumber\\[1ex]
&&\hspace{1cm}=\;\;g^k_b\vec{\ell}^{\,ij}\wedge\nabla_a\bp-g^k_a\vec{\ell}^{\,ij}\wedge\nabla_b\bp+(\vec{\ell}^{\,ij}\cdot\nabla_a\bp)\vec{\eta}^{\,k}_{\,\,\;b}-(\vec{\ell}^{\,ij}\cdot\nabla_b\bp)\vec{\eta}^{\,k}_{\,\;\,a}\:.
\end{eqnarray}
Accordingly, we find on one hand:
\begin{eqnarray}\label{ahah2}
\bD\stackrel{\!\!\!\bullet}{\res}\vec{\eta}&=&2\big(\vec{\ell}\stackrel{\wedge}{\wedge}d\bp\big)\stackrel{\!\!\!\bullet}{\res}\vec{\eta}\nonumber\\
&=&\big(\vec{\ell}^{\,ab}\wedge\nabla^c\bp+\vec{\ell}^{\,ca}\wedge\nabla^b\bp +\vec{\ell}^{\,bc}\wedge\nabla^a\bp \big)\bullet(\nabla_a\bp\wedge\nabla_b\bp)\nonumber\\[1ex]
&=&\big(\vec{\ell}^{\,ab}\wedge\nabla^c\bp+2\vec{\ell}^{\,ca}\wedge\nabla^b\bp \big)\bullet(\nabla_a\bp\wedge\nabla_b\bp)\nonumber\\[1ex]
&=&2\vec{\ell}^{\,ac}\wedge\nabla_a\bp+2(\vec{\ell}^{\,ab}\cdot\nabla_a\bp)\vec{\eta}^{\,c}_{\:\:\:b}+6\vec{\ell}^{\,ca}\wedge\nabla_a\bp-2(\vec{\ell}^{\,ca}\cdot\nabla_b\bp)\vec{\eta}^{\,b}_{\:\:\:a}\nonumber\\[1ex]
&=&-\,4\bCC-2\vec{\eta}\res A+2(\vec{\ell}^{\,ca}\cdot\nabla^b\bp)\vec{\eta}_{ab}\:.
\end{eqnarray}
On the other hand, we have
\begin{eqnarray}\label{ahah3}
\vec{\eta}\stackrel{\!\!\!\bullet}{\res}\bCC&=&\vec{\eta}\stackrel{\!\!\!\bullet}{\res}\big(\vec{\,\ell}\stackrel{\!\!\!\wedge}{\res}d\bp\big)\;\;=\;\;\vec{\eta}^{\,bc}\bullet\big(\vec{\ell}_{ab}\wedge\nabla^a\bp\big)\nonumber\\[1ex]
&=&(\nabla^b\bp\wedge\nabla^c\bp)\bullet\big(\vec{\ell}_{ab}\wedge\nabla^a\bp\big)\nonumber\\[1ex]
&=&\nabla_b\bp\wedge\vec{\ell}^{\,cb}+(\vec{\ell}^{\,ab}\cdot\nabla^c\bp)\vec{\eta}_{ab}+\big(\vec{\ell}_{ab}\cdot\nabla^b\bp\big)\vec{\eta}^{\,ca}\nonumber\\[1ex]
&=&\bCC+\vec{\eta}\res A+\big(\vec{\ell}^{\,ab}\cdot\nabla^c\bp\big)\vec{\eta}_{ab}\:.
\end{eqnarray}
Finally, we compute
\begin{eqnarray}\label{ahah4}
B\res\vec{\eta}&=&2\big(\vec{\ell}\stackrel{\cdot}{\wedge}d\bp\big)\res\vec{\eta}\nonumber\\
&=&\big(\vec{\ell}^{\,ab}\cdot\nabla^c\bp+\vec{\ell}^{\,ca}\cdot\nabla^b\bp +\vec{\ell}^{\,bc}\cdot\nabla^a\bp \big)\vec{\eta}_{ab}\nonumber\\[1ex]
&=&2(\vec{\ell}^{\,ca}\cdot\nabla^b\bp)\vec{\eta}_{ab}+(\vec{\ell}^{\,ab}\cdot\nabla^c\bp \big)\vec{\eta}_{ab}\:.
\end{eqnarray}
Combining together (\ref{ahah2})-(\ref{ahah4}) yields now the desired (\ref{clarisse}):
\be\label{clarissee}
 -\,3\bCC\;=\;\vec{\eta}\stackrel{\!\!\!\bullet}{\res}\bCC+\bD\stackrel{\!\!\!\bullet}{\res}\vec{\eta}+\vec{\eta}\res A-B\res\vec{\eta}\:.
 \ee

To prove (\ref{clarisse2}), we proceed analogously:
\begin{eqnarray*}
\vec{\eta}\stackrel{\!\!\!\cdot}{\res}\bCC&=&\vec{\eta}^{\,bc}\cdot(\vec{\ell}_{ab}\wedge\nabla^a\bp)\;\;=\;\;(\nabla^b\bp\wedge\nabla^c\bp)\cdot(\vec{\ell}_{ab}\wedge\nabla^a\bp)\\[1ex]
&=&\vec{\ell}_{cb}\cdot\nabla^b\bp\\
&=&-A
\end{eqnarray*}
and
\begin{eqnarray*}
\bD\stackrel{\!\!\!\cdot}{\res}\vec{\eta}&=&2\big(\vec{\ell}\stackrel{\wedge}{\wedge}d\bp\big)\stackrel{\!\!\!\cdot}{\res}\vec{\eta}\nonumber\\
&=&\big(\vec{\ell}^{\,ab}\wedge\nabla^c\bp+\vec{\ell}^{\,ca}\wedge\nabla^b\bp +\vec{\ell}^{\,bc}\wedge\nabla^a\bp \big)\cdot(\nabla_a\bp\wedge\nabla_b\bp)\nonumber\\[1ex]
&=&\big(\vec{\ell}^{\,ab}\wedge\nabla^c\bp+2\vec{\ell}^{\,ca}\wedge\nabla^b\bp \big)\cdot(\nabla_a\bp\wedge\nabla_b\bp)\nonumber\\[1ex]
&=&-4\vec{\ell}^{\,ac}\cdot\nabla_a\bp\nonumber\\[1ex]
&=&-4A\:.
\end{eqnarray*}
Hence we recover the announced (\ref{clarisse2}):
\bes\label{guy1}
3A\;=\;\vec{\eta}\stackrel{\!\!\!\cdot}{\res}\bCC-\bD\stackrel{\!\!\!\cdot}{\res}\vec{\eta}\:.
\ees
%
\hfill$\blacksquare$\\

We apply this proposition with the choice $\vec{\ell}=\bL$ and
$$
A\;=\;d^\star S\quad,\quad B\;=\;dS\quad,\quad\bCC\;=\;d^\star\bR-\bU\quad,\quad\bD\;=\;d\bR\:.
$$
Doing so yields the system
\begin{equation}\label{sysSR1}
\left\{\begin{array}{rcl}
-3d^\star\bR&=&\vec{\eta}\stackrel{\!\!\!\bullet}{\res}d^\star\bR+d\bR\stackrel{\!\!\!\bullet}{\res}\vec{\eta}+\vec{\eta}\res d^\star S-dS\res\vec{\eta}-\vec{Y}\\
3d^\star S&=&\vec{\eta}\stackrel{\!\!\!\cdot}{\res}d^\star\bR-d\bR\stackrel{\!\!\!\cdot}{\res}\vec{\eta}-Z\:,
\end{array}\right. 
\end{equation}
where for notational convenience we have set
$$
\vec{Y}\;:=\;3\bU+\vec{\eta}\stackrel{\!\!\!\bullet}{\res}\bU\qquad\text{and}\qquad Z\;:=\;\vec{\eta}\stackrel{\!\!\!\cdot}{\res}\bU\:.
$$
To further proceed, it is necessary to consider the function $\bF$ carefully. According to the results from Section~\ref{Noether}, we see that 
$$
\bF^{ij}\;=\;-\dfrac{1}{2}g^{ij}\Delta_\perp\bH+\mathcal{O}(|\bh|^3)\:,
$$
so that, from (\ref{defU}), 
$$
\bU^j\;=\;\dfrac{1}{2}\Delta_\perp\bH\wedge\nabla^j\bp+\mathcal{O}\big(|\bh|^3+|\bH||\pro d\bH|\big)\:.
$$
This implies easily that
$$
Z\;=\;\mathcal{O}\big(|\bh|^3+|\bH||\pro d\bH|\big)\:.
$$
Moreover, as $\vec{v}:=\frac{1}{2}\Delta_\perp\bH$ is a normal-valued scalar function, it holds
\begin{eqnarray*}
3\bv\wedge\nabla^j\bp+\vec{\eta}^{ij}\bullet(\bv\wedge\nabla_i\bp)&=&3\bv\wedge\nabla^j\bp+(\nabla^i\bp\wedge\nabla^j\bp)\bullet(\vec{v}\wedge\nabla_i\bp)\\
&=&3\bv\wedge\nabla^j\bp+\big(\nabla^j\bp\wedge\bv-4\nabla^j\bp\wedge\bv)\\
&=&6\vec{v}\wedge\nabla^j\bp\:.
\end{eqnarray*}
This implies
$$
\vec{Y}^j\;=\;3\Delta_\perp\bH\wedge\nabla^j\bp+\mathcal{O}\big(|\bh|^3+|\bH||\pro d\bH|\big)\:.
$$
We now call again upon the 1-vector-valued 1-form 
$$
\bX^i\;=\;\nabla^i\bH-2(|\bH|^2g^{ij}-\bH\cdot\bh^{ij})\nabla_j\bp\:.
$$
We have seen in (\ref{yab}) that $d^{\star}\bX=\Delta_\perp\bH+\mathcal{O}(|\bh|^3)$. Hence
\begin{eqnarray*}
d^\star(\bX\stackrel{\wedge}{\wedge}d\bp)&=&\nabla_i\big(\bX^i\wedge\nabla^j\bp-\bX^j\wedge\nabla^i\bp)\\
&=&d^\star\bX\wedge d\bp-\nabla_i\bX^j\wedge\nabla^i\bp+\mathcal{O}\big(|\bh||\bX|\big)\\
&=&d^\star\bX\wedge d\bp-(\nabla^i\bX^j-\nabla^j\bX^i)\wedge\nabla_i\bp-\nabla_j(\bX^i\wedge\nabla_i\bp)+\mathcal{O}\big(|\bh||\bX|\big)\\
&=&d^\star\bX\wedge d\bp-d\langle\bX\stackrel{\wedge}{,}d\bp\rangle+\mathcal{O}\big(|d\bX|+|\bh||\bX|\big)\:.
\end{eqnarray*}
Thus we can express
$$
\vec{Y}\;=\;3d^\star(\bX\stackrel{\wedge}{\wedge}d\bp)+3d\langle\bX\stackrel{\wedge}{,}d\bp\rangle+\mathcal{O}\big(|\bh|^3+|\bh||\pro d\bH|\big)\:,
$$
where we have used that $d\bX=\mathcal{O}(|\bh|^3+|\bh||\pro d\bH|)$, which easily follows from Codazzi. \\

This information may now be imported into (\ref{sysSR1}) so as to write
\begin{equation}\label{sysSR2}
\left\{\begin{array}{rcl}
-3d^\star(\bR-\bX\stackrel{\wedge}{\wedge}d\bp)+3d\langle\bX\stackrel{\wedge}{,}d\bp\rangle&=&\vec{\eta}\stackrel{\!\!\!\bullet}{\res}d^\star\bR+d\bR\stackrel{\!\!\!\bullet}{\res}\vec{\eta}+\vec{\eta}\res d^\star S-dS\res\vec{\eta}+\mathcal{O}\big(|\bh|^3+|\bh||\pro d\bH|\big)\\
3d^\star S&=&\vec{\eta}\stackrel{\!\!\!\cdot}{\res}d^\star\bR-d\bR\stackrel{\!\!\!\cdot}{\res}\vec{\eta}+\mathcal{O}\big(|\bh|^3+|\bH||\pro d\bH|\big)\:,
\end{array}\right. 
\end{equation}
which completes the proof of Theorem~\ref{ThSR}. \\

To complete the proof of Corollary~\ref{COSR}, we call upon the following proposition:
\begin{Prop}\label{strucRS}
Let $A$ and $B$ be two real-valued two-forms. Then
$$
A\res d^\star B-dB\res A\;=\;d^\star(A\odot B)-d\langle A,B\rangle+\mathcal{O}(|\nabla A||B|)\:.
$$
Similarly, if $\bA$ and $\bB$ are two two-forms with values in $\Lambda^2(\R^{m\ge5})$, then
$$
\bA\stackrel{\!\!\!\cdot}{\res} d^\star\bB-d\bB\stackrel{\!\!\!\cdot}{\res}\bA\;=\;d^\star(\bA\stackrel{\cdot}{\odot}\bB)-d\langle\bA\stackrel{\cdot}{,}\bB\rangle+\mathcal{O}(|\nabla\bA||\bB|)\:.
$$
and
$$
\bA\stackrel{\!\!\!\bul}{\res} d^\star\bB+d\bB\stackrel{\!\!\!\bul}{\res}\bA\;=\;d^\star(\bA\stackrel{\bul}{\odot}\bB)-d\langle\bA\stackrel{\bul}{,}\bB\rangle+\mathcal{O}(|\nabla\bA||\bB|)\:.
$$

\end{Prop}
{\bf Proof.} Suppose first that $A$ and $B$ are real-valued 2-forms. By direct computation, we have\footnote{I use the symbol ${\mathlarger{\mathlarger{\mathlarger{\lrcorner}}}}$ to distinguish between the left and right actions of contraction. This is of importance when we consider 2-forms with values in a noncommutative space (that of 2-vectors for instance).}
\begin{eqnarray*}
2A_{ai}\nabla^aB^{ib}&=&A_{ai}\nabla^aB^{ib}+A_{ai}\nabla^iB^{ba}\;\;=\;\;A_{ai}(dB)^{aib}-A_{ai}\nabla^{b}B^{ai}\\[1ex]
&=&2A\mathlarger{\mathlarger{\mathlarger{\mathlarger{\lrcorner}}}}dB-\nabla^b(A_{ai}B^{ai})+\mathcal{O}(|\nabla A||B|)\\[1ex]
&=&2A\mathlarger{\mathlarger{\mathlarger{\mathlarger{\lrcorner}}}}dB-2d\langle A,B\rangle+\mathcal{O}(|\nabla A||B|)\:.
\end{eqnarray*}
Hence
\begin{eqnarray}\label{choisi}
-d^\star(A\odot B)&\equiv&\nabla^a\big(A_{ai}B^{ib}-A^{ib}B_{ai}\big)\nonumber\\[1ex]
&=&-A^{ib}\nabla^aB_{ai}+A_{ai}\nabla^aB^{ib}   +\mathcal{O}(|\nabla A||B|)\nonumber\\[1ex]
&=&-\,A\res d^\star B+A\mathlarger{\mathlarger{\mathlarger{\mathlarger{\lrcorner}}}}dB-d\langle A,B\rangle+\mathcal{O}(|\nabla A||B|)\:,
\end{eqnarray}
which yields the first identity. \\
Suppose now that $A$ and $B$ map into the space of 2-vectors of $\R^m$. To be consistent with our notation, we render them as $\bA$ and $\bB$. Applying first in (\ref{choisi}) the operation $\cdot$ to the vectorial-components of the elements involved yields
$$
d^\star(\bA\stackrel{\cdot}{\odot}\bB)-d\langle\bA\stackrel{\cdot}{,}\bB\rangle\;=\;\bA\stackrel{\!\!\!\cdot}{\res} d^\star\bB-d\bB\stackrel{\!\!\!\cdot}{\res}\bA+\mathcal{O}(|\nabla\bA||\bB|)\:,
$$
where we have used that $\cdot$ commutes on 2-vectors. On the other hand, applying in (\ref{choisi}) the operation $\bullet$ to the vectorial-components of the elements involved yields
$$
d^\star(\bA\stackrel{\bul}{\odot}\bB)-d\langle\bA\stackrel{\bul}{,}\bB\rangle\;=\;\bA\stackrel{\!\!\!\bul}{\res} d^\star\bB+d\bB\stackrel{\!\!\!\bul}{\res}\bA+\mathcal{O}(|\nabla\bA||\bB|)\:,
$$
where we have used that $\bul$ anticommutes on 2-vectors.

\hfill$\blacksquare$\\

We may now apply this proposition directly to the system (\ref{sysSR2}) with $\bA=\vec{\eta}$:
\begin{eqnarray*}
&&d^\star\big(3\bR-3\bX\stackrel{\wedge}{\wedge}d\bp+\vec{\eta}\stackrel{\bullet}{\odot}\bR+\vec{\eta}\odot S\big)-d\big(\langle\vec{\eta}\stackrel{\bullet}{,}\bR\rangle+\langle\vec{\eta},S\rangle+3\langle\bX\stackrel{\wedge}{,}d\bp\rangle\big)\\
&&\hspace{5cm}=\:\mathcal{O}\big(|\bh|^3+|\bh||\pro d\bH|+|\nabla\vec{\eta}||\bR|+|\nabla\vec{\eta}||S|\big)\:,
\end{eqnarray*}
as well as
$$
d^\star\big(3S-\vec{\eta}\stackrel{\cdot}{\odot}\bR\big)+d\langle\vec{\eta}\stackrel{\cdot}{,}\bR\rangle\;=\;\mathcal{O}\big(|\bh|^3+|\bH||\pro d\bH|+|\nabla\vec{\eta}||\bR|\big)\:, 
$$
Since $|\nabla\vec{\eta}|=\mathcal{O}(|\bh|)$, we reach the statement of Corollary~\ref{COSR} by applying $d^{\star}$ to both sides of the first identity.

\subsection{Return to geometry: proof of Theorem~\ref{Threturn}}

As in the previous section, the proof relies on algebraic identities. 

\begin{Prop}\label{last}
Let $\vec{R}\in\Lambda^2(\Lambda^2)$ and $S\in\Lambda^2(\Lambda^0)$. 
Denoting by $\pro$ the projection onto the normal subspace, it holds
$$
\pro\, \Big(d^\star\vec{R}\stackrel{\res}{\res}d\bp\Big)\;\equiv\;-\pro\, d^\star\Big(\vec{R}\stackrel{\res}{\res}d\bp\Big)\;=\;-\pro\,d^\star\Big(\langle\vec{\eta}\stackrel{\bullet}{,}\vec{R}\rangle\res d\bp\Big)+\mathcal{O}(|\vec{R}||\bh|)\:.
$$
Furthermore, we have
$$
\pro\,d^\star\big(\langle\vec{\eta}\,,S\rangle\res d\bp\big)\;=\;-\pro\big(d^\star S\res d\bp\big)\;=\;\vec{0}\:.
$$
\end{Prop}
{\bf Proof.} We decompose 
%
%
%
%
%
%
%
%
%
$\vec{R}=\vec{R}_{\bn}+\vec{R}_{\vec{t}}$, where $\vec{R}_{\vec{t}}\,$ involves only linear combinations of 2-vectors of the type $\nabla^a\bp\wedge\nabla^b\bp$. In particular, it follows that 
$$
\pro\Big(\vec{R}_{\vec{t}}\stackrel{\res}{\res}d\bp\Big)\;=\;\vec{0}\:,
$$
thereby yielding that
$$
\vec{R}_{\vec{t}}\stackrel{\res}{\res}d\bp\;=\;u_{ai}\nabla^i\bp\,dx^a\:,
$$
for some coefficients $u_{ai}$ {\it which we do not assume to be either symmetric or antisymmetric}. Hence
$$
\pro\,d^\star\Big(\vec{R}_{\vec{t}}\stackrel{\res}{\res}d\bp\Big)\;=\;\mathcal{O}(|u||\bh|)\:.
$$
As 
$$
|u|\;=\;\mathcal{O}(|\vec{R}_{\vec{t}}|)\;=\;\mathcal{O}(|\vec{R}|)\:,
$$
we arrive at
\be\label{lolo1}
\pro\, d^\star\Big(\vec{R}\stackrel{\res}{\res}d\bp\Big)\;=\;\pro\, d^\star\Big(\vec{R}_{\bn}\stackrel{\res}{\res}d\bp\Big)+\mathcal{O}(|\vec{R}||\bh|)\:.
\ee
\medskip
By the same token, we note that $\langle\vec{\eta}\stackrel{\bullet}{,}\vec{R}_{\vec{t}}\rangle$ is also a linear combination of 2-vectors of the type $\nabla^a\bp\wedge\nabla^b\bp$. Just as above, we obtain that
\be\label{lolo2}
\pro\,d^\star\Big(\langle\vec{\eta}\stackrel{\bullet}{,}\vec{R}_{\vec{t}}\rangle\res d\bp\Big)\;=\;\mathcal{O}(|\vec{R}||\bh|)\:.
\ee
\medskip
We next project $\vec{R}_{\bn}$ to write
$$
\big(\vec{R}_{\bn}\big)_{ab}\;=\;[\vec{f}_{ab}]_j\,\wedge\nabla^j\bp+\vec{C}\:,
$$
where for each fixed index $j$, the coefficient $[\vec{f}_{ab}]_j$ is that of a normal-valued two-form; and $\vec{C}$ is a 2-form spanned by the wedge product of normal-valued 1-forms. By direct computation, we check on one hand that
\begin{eqnarray*}
\langle\vec{\eta}\stackrel{\bullet}{,}\vec{R}_{\bn}\rangle&\equiv&\dfrac{1}{2}\vec{\eta}^{\,ab}\bullet\big(\vec{R}_{\bn}\big)_{ab}\;\;=\;\;\dfrac{1}{2}(\nabla^a\bp\wedge\nabla^b\bp)\bullet ([\vec{f}_{ab}]_j\wedge\nabla^j\bp)\\
&=&[\vec{f}_{ab}]^a\wedge\nabla^b\bp\:,
\end{eqnarray*}
so that
$$
\langle\vec{\eta}\stackrel{\bullet}{,}\vec{R}_{\bn}\rangle\res \nabla_j\bp\;=\;-[\vec{f}_{aj}]^a\:;
$$
and on the other hand that
\begin{eqnarray*}
\big(\vec{R}_{\bn}\stackrel{\res}{\res}d\bp\big)_j&=&([\vec{f}_{aj}]_b\wedge\nabla^b\bp)\res\nabla^a\bp\;\;=\;\;-[\vec{f}_{aj}]^a\:.
\end{eqnarray*}
Accordingly,
\bes
\vec{R}_{\bn}\stackrel{\res}{\res}d\bp\;=\;\langle\vec{\eta}\stackrel{\bullet}{,}\vec{R}_{\bn}\rangle\res d\bp\:.
\ees
It follows from the latter and (\ref{lolo1}), (\ref{lolo2}) that
$$
\pro\, d^\star\Big(\vec{R}\stackrel{\res}{\res}d\bp\Big)\;=\;\pro\,d^\star\Big(\langle\vec{\eta}\stackrel{\bullet}{,}\vec{R}\rangle\res d\bp\Big)+\mathcal{O}(|\vec{R}||\bh|)\:,
$$
as claimed.
\medskip
To finish the proof, it remains to show the second identity. This is also done by direct computation:
\begin{eqnarray*}
\langle\vec{\eta}\,,S\rangle\res d\bp&=&\dfrac{1}{2}S_{ab}(\nabla^a\bp\wedge\nabla^b\bp)\res\nabla^i\bp\;=\;S^{ai}\,\nabla_a\bp\;=\;S\res d\bp\:,
\end{eqnarray*}
which immediately implies the announced statement. 

$\hfill\blacksquare$\\

By definition, (\ref{defR}) yields
\begin{eqnarray}\label{flip22}
d^\star\bR\stackrel{\res}{\res}d\bp&=&(\bL\stackrel{\!\!\!\wedge}{\res}d\bp)\stackrel{\res}{\res}d\bp+\bU\stackrel{\res}{\res}d\bp\nonumber\\
&=&(\bL\stackrel{\!\!\!\wedge}{\res}d\bp)\stackrel{\res}{\res}d\bp-2d^\star\bX\:,
\end{eqnarray}
where we have used (\ref{huma4}). \\
We note next that
\begin{eqnarray*}
(\bL\stackrel{\!\!\!\wedge}{\res}d\bp)\stackrel{\res}{\res}d\bp+(\bL\stackrel{\!\!\!\cdot}{\res}d\bp)\stackrel{}{\res}d\bp&=&\big(\bL^{ab}\wedge\nabla_a\bp\big)\res\nabla_b\bp+\big(\bL^{ab}\cdot\nabla_a\bp\big)\nabla_b\bp\nonumber\\
&=&\big(\bL^{ab}\cdot\nabla_b\bp\big)\nabla_a\bp+\big(\bL^{ab}\cdot\nabla_a\bp\big)\nabla_b\bp\nonumber\\
&=&\vec{0}\:.
\end{eqnarray*}
This shows that
\bes
(\bL\stackrel{\!\!\!\wedge}{\res}d\bp)\stackrel{\res}{\res}d\bp\;=\;-d^\star S\res d\bp\:.
\ees
Hence now (\ref{flip22}) becomes
\bes
d^\star\bR\stackrel{\res}{\res}d\bp+d^\star S\res d\bp\;=\;-2d^\star\bX\:.
\ees
Note that this is indeed an equation in divergence form. Combining it to Proposition~\ref{last} gives finally:
\be\label{kus}
2d^\star\bX\;=\;\pro\,d^\star\Big(\big(\langle\vec{\eta}\stackrel{\bullet}{,}\vec{R}\rangle+\langle \vec{\eta},S\rangle\big)\res d\bp\Big)+\mathcal{O}\big(|\vec{R}||\bh|\big)\:,
\ee
where we have used that $d^\star\bX$ is normal-valued, as explained in Lemma~\ref{lem4}. \\

Recall that
$$
\bX^i\;=\;\pro\nabla^i\bH-\big(2|\bH|^2g^{ij}-\bH\cdot\bh^{ij}\big)\nabla_j\bp\
$$
satisfies $d^\star\bX=\Delta_\perp\bH+\mathcal{O}(|\bh|^3)$, so that (\ref{kus}) yields
\bes
2\Delta_\perp\bH\;=\;\pro\,d^\star\Big(\big(\langle\vec{\eta}\stackrel{\bullet}{,}\vec{R}\rangle+\langle \vec{\eta},S\rangle\big)\res d\bp\Big)+\mathcal{O}\big(|\vec{R}||\bh|+|\bh|^3\big)\:.
\ees
Moreover,
$$
\langle\bX\stackrel{\wedge}{,}d\bp\rangle\res d\bp\;=\;(\pro\nabla^i\bH\wedge\nabla_i\bp)\res\nabla^j\bp\;=\;-\pro\nabla^j\bH\:,
$$
and thus 
\bes
-\Delta_\perp\bH\;=\;\pro\,d^\star\Big(\big(\langle\vec{\eta}\stackrel{\bullet}{,}\vec{R}\rangle+\langle \vec{\eta},S\rangle+3\langle\bX\stackrel{\wedge}{,}d\bp\rangle\big)\res d\bp\Big)+\mathcal{O}\big(|\vec{R}||\bh|+|\bh|^3\big)\:,
\ees
which is the content of Theorem~\ref{Threturn}.

\bigskip

\setcounter{equation}{0} 
\reset

\renewcommand{\theequation}{A.\arabic{equation}}
\renewcommand{\theTh}{A.\arabic{Th}}
\renewcommand{\theProp}{A.\arabic{Prop}}
\renewcommand{\theLm}{A.\arabic{Lm}}
\renewcommand{\theCo}{A.\arabic{Co}}
\renewcommand{\theRm}{A.\arabic{Rm}}
\renewcommand{\theequation}{A.\arabic{equation}}
\setcounter{equation}{0} 
\reset
\appendix
\section{Appendix}

\subsection{Notational Conventions}\label{nota}
\reset

We fix a basis $\{e^i\}_{i=1,\cdots,N}$ of $\mathbb{R}^N$ and an orientation. Suppose $\mathbb{R}^N$ is equipped with a metric $\frak{g}$. Let $\Lambda^k(\R^N)$ be the corresponding vector space of $k$-forms. We will need component expansions of certain operations involving forms. A typical element $A\in\Lambda^p$ may be expanded in tensorial notation as
$$
A\;=\;\dfrac{1}{p!}\sum_{\text{all indices}}A_{i_1\ldots i_p}\,e^{i_1}\wedge\ldots\wedge e^{i_p}\:,
$$
where $A_{i_1\ldots i_p}$ is totally antisymmetric. For notational convenience, we will render this as
$$
A\;\equiv\;A_{i_1\ldots i_p}\:.
$$
Forms are tensors: we raise and lower indices as usual using the metric $\frak{g}$. Of course, as forms are anti-symmetric, the relative position of indices is essential: we follow standard conventions. \\
The inner product on forms will be denoted
$$
\langle A,B\rangle_\frak{g}\;=\;\dfrac{1}{p!}A^{i_1\ldots i_p}B_{i_1\ldots i_p}\qquad\text{for}\quad A\in\Lambda^p\:\:,\:\: B\in\Lambda^p\:.
$$
The bracket symbol $[\ldots]$ will denote \underline{scaled} total permutations. We conventionally set the wedge product to be
$$
A\wedge B\;\equiv\;\dfrac{(p+q)!}{p!q!}A_{[i_1\ldots i_p}B_{j_1\ldots j_q]}\qquad\text{for}\quad A\in\Lambda^p\:\:,\:\: B\in\Lambda^q\:.
$$
The {\it interior multiplication} (restricted to our interest) is the bilinear operation $\res\;:\;\Lambda ^p(\mathbb{R}^N)\times \Lambda^{p-1}(\mathbb{R}^N)\rightarrow \Lambda^{1}(\mathbb{R}^N)$ which satisfies in particular
$$
A\res B\;\equiv\;\dfrac{1}{(p-1)!}A^{i_1\ldots i_p}B_{i_1\ldots i_{p-1}}\qquad\text{for}\quad A\in\Lambda^p\:\:,\:\: B\in\Lambda^{p-1}\:.
$$
\noindent
The {\it first-order contraction operator}   $\bullet: \Lambda^2(\mathbb{R}^N)\times \Lambda^q(\mathbb{R}^N)\rightarrow \Lambda^{q}(\mathbb{R}^N)$ is defined as follows. For a $2$-form ${A}$ and a $1$-form ${B}$, we set 
$${A}\bullet {B}:={A}\res B\:.$$
For a $p$-form ${B}$ and a $q$-form ${C}$, we have
$${A}\bullet ({B}\wedge{C}):=({A}\bullet {B})\wedge {C}+ (-1)^{pq}({A}\bullet{C})\wedge{B}\:.$$
This operation will be used for forms with values in $\Lambda^p(\R^{m\ge5})$.\\
We will also need an operation taking a pair of two-forms to a two-form, namely:
$$
A\odot B\;\equiv\;A^{i_2j}B_{j}^{\:\:i_1}-A^{i_1j}B_{j}^{\:\:i_2}\qquad\text{for}\quad A, B\in\Lambda^2\:.
$$
\noindent
Being given an orientation and a metric, we can equip our manifold with the volume form $\epsilon\in\Lambda^N$ via $\epsilon^{i_1\ldots i_N}:=|\frak{g}|^{-1/2}\mbox{sign} \begin{pmatrix}
1 & 2 & \cdots & N\\
i_1 &i_2 & \cdots & i_N
\end{pmatrix}$. With it, we define the Hodge star operator $\star_\frak{g}:\Lambda^p\rightarrow\Lambda^{N-p}$ as usual. In components, we have
$$
\star_\frak{g}A\;\equiv\;\dfrac{1}{p!}\,\epsilon^{i_1\ldots i_{N-p}j_1\ldots j_p}A_{j_1\ldots j_p}\qquad\text{for}\quad A\in\Lambda^p\:.
$$
Note that
$$
\star_\frak{g}\star_\frak{g} A=(-1)^{p(N-p)}A\:.
$$

\noindent
If a connection $\nabla$ compatible with the metric $\frak{g}$ is provided, we define in accordance with the Leibniz rule, the covariant exterior derivative $d$ of a $p$-form $A$ in components as
$$dA\;:=\;(p+1)\nabla_{[i_0}A_{i_1...i_p]}\;\in\;\Lambda^{p+1}\:.$$
We also define the codifferential $d^{\star_\frak{g}}:=\star_\frak{g} d\star_\frak{g}$, whose components take the expression
$$
d^{\star_\frak{g}} A=\nabla^{j}A_{ji_1...i_{p-1}}\;\in\;\Lambda^{p-1}\:.
$$
For any $p$-form $A$, it holds
$$
d^2A\;=\;0\;=\;(d^{\star_\frak{g}})^2A\:.
$$

\medskip

\noindent
We will work within the context of multivector-valued forms. Our objects will be elements of $\Lambda^p(\R^4,\Lambda^q(\R^5))$. To distinguish operations taking place in the parameter space $(\R^4,g)$ from operations taking place in the ambient space $(\R^{m},g_{\text{Eucl}})$, we will superimpose the operators defined above with the convention that the one on the bottom acts on parameter space, while the one at the top acts on ambient space. 
Typically, consider the two-vector-valued two-form $\vec{\eta}\in\Lambda^2(\R^4,\Lambda^2(\R^m))$ whose components are given by
$$
\vec{\eta}_{ij}\;=\;\nabla_i\bp\wedge\nabla_j\bp\:.
$$
We will write more succinctly
$$
\vec{\eta}\;=:\;\dfrac{1}{2}d\bp\stackrel{\wedge}{\wedge}d\bp\:.
$$
Similar quantities are defined following the same logic. Naturally, for multivector 0-forms, the bottom symbol will be omitted, while for 0-vector-valued forms, the symbol on the top will be omitted instead. The notation $\langle\cdot\,,\cdot\rangle$ is reserved for $\langle\cdot\,,\cdot\rangle_g$. The symbol $\star$ is reserved for $\star_g$. In other words, the Hodge star operator used in the present article is intrinsic. 
Similarly, $\Delta$ will stand for $\Delta $, the Laplace-Beltrami operator (or Bochner Laplacian when acting on the components of tensors). 

\subsection{Noether Fields}\label{Noether}
\reset

\subsubsection{The energy $\int|\pro d\bH|^2$}\label{vardh}

We consider a variation of the type
\be\label{varo}
\delta\bp\;=\;A^s\nabla_s\bp+\bB\:,
\ee
where $\bB$ is a normal vector. \\
The metric varies according to:
\begin{equation}\label{varmet}
\delta g^{ij}\;=\;-\nabla^iA^j-\nabla^jA^i+2\bB\cdot\bh^{ij}\qquad\text{and}\qquad\delta|g|^{1/2}\;=\;|g|^{1/2}\big(\nabla_jA^j-4\bB\cdot\bH\big)\:.
\end{equation}
Note that
\begin{equation*}
\delta\nabla_j\bp\;=\;\nabla_j\delta\bp\;=\;\nabla_j\big(A^s\nabla_s\bp+\bB\big)\;=\;A^s\bh_{sj}+\nabla_jA^s\nabla_s\bp+\nabla_j\bB\:.
\end{equation*}
Hence, using Codazzi,
\begin{eqnarray}\label{eq0}
\pro\nabla^j\delta\nabla_j\bp&=&\bh_{sj}\nabla^jA^s+A^s\pro\nabla^j\bh_{sj}+\bh^j_s\nabla_jA^s+\pro\Delta\bB\nonumber\\[1ex]
&=&2\bh_{sj}\nabla^jA^s+4A^s\pro\nabla_s\bH+\pro\Delta\bB\:.
\end{eqnarray}
We next 
compute the normal part of the variation of $\bH$, using (\ref{varmet}) and (\ref{eq0}):
\begin{eqnarray}\label{pindelh}
\pro\delta\bH&=&\dfrac{1}{4}\pro\delta\nabla^j\nabla_j\bp\;\;=\;\;\dfrac{1}{4}\pro\delta\big(g^{ij}\nabla_{i}\nabla_j\bp\big)\nonumber\\[1ex]
&=&\dfrac{1}{4}\big(\bh_{ij}\delta g^{ij}+\pro\nabla^j\delta\nabla_j\bp\big)\nonumber\\[1ex]
&=&\dfrac{1}{4}\bh_{ij}\big(-\nabla^iA^j-\nabla^jA^i+2\bB\cdot\bh^{ij}\big)+\dfrac{1}{4}\big(2\bh_{sj}\nabla^jA^s+4A^s\pro\nabla_s\bH+\pro\Delta\bB\big)\nonumber\\[1ex]
&=&\dfrac{1}{4}\big(4A^s\pro\nabla_s\bH+2(\bB\cdot\bh^{ij})\bh_{ij}+\pro\Delta\bB\big)\nonumber\\[1ex]
&=&\dfrac{1}{4}\big(4A^s\pro\nabla_s\bH+(\bB\cdot\bh^{ij})\bh_{ij}+\Delta_\perp\bB\big)\:.
\end{eqnarray}

Given a tensor field $\bT_j$, the variation of $|\bT|^2$ is found from the variation of $\bT_j$ as follows:
\begin{eqnarray*}
\delta|\bT|^2&=&\delta(g^{ij}\bT_i\cdot\bT_j)\;\;=\;\;(\bT_i\cdot\bT_j)\delta g^{ij}+2\bT^j\cdot\delta\bT_j\nonumber\\[1ex]
&=&2\bT_i\cdot\bT_j\big(-\nabla^iA^j+\bB\cdot\bh^{ij}\big)+2\bT^j\cdot\delta\bT_j\:,
\end{eqnarray*}
where we have used (\ref{varmet}). Using again (\ref{varmet}), we have
\begin{equation}\label{vary}
|g|^{-1/2}\delta\big(|\bT|^2|g|^{1/2}\big)\;=\;2\bT_i\cdot\bT_j\big(-\nabla^iA^j+\bB\cdot\bh^{ij}\big)+2\bT^j\cdot\delta\bT_j+|\bT|^2\big(\nabla_jA^j-4\bB\cdot\bH\big)\:.
\end{equation}

\eject

We will now vary $|\bT|^2|g|^{1/2}$ for the choice $\bT_j=\pro\nabla_j\bH$. We first note that
\begin{eqnarray}\label{eq10}
&&\pro\nabla^j\bH\cdot\delta\pro\nabla_j\bH\;\;=\;\;\pro\nabla^j\bH\cdot\delta\nabla_j\bH-\pro\nabla^j\bH\cdot\delta\pi_T\nabla_j\bH\nonumber\\[1ex]
&=&\pro\nabla^j\bH\cdot\nabla_j\delta\bH-\pro\nabla^j\bH\cdot\delta\pi_T\nabla_j\bH\nonumber\\[1ex]
&=&\pro\nabla^j\bH\cdot\nabla_j\pro\delta\bH+\pro\nabla^j\bH\cdot\nabla_j\pi_T\delta\bH-\pro\nabla^j\bH\cdot\delta\pi_T\nabla_j\bH\nonumber\\[1ex]
&=&\pro\nabla^j\bH\cdot\nabla_j\pro\delta\bH+\pro\nabla^j\bH\cdot\Big[\pro\nabla_j\big[(\nabla_i\bp\cdot\delta\bH)\nabla^i\bp\big]-\delta\big[(\nabla^i\bp\cdot\nabla_j\bH)\nabla_i\bp   \big]  \Big]\nonumber\\[1ex]
&=&\pro\nabla^j\bH\cdot\nabla_j\pro\delta\bH+\nabla^j\bH\cdot\Big[-\pro\nabla_j\big[(\bH\cdot\nabla_i\delta\bp)\nabla^i\bp\big]+\pro\delta\big[(\bH\cdot\bh^i_j)\nabla_i\bp   \big]  \Big]\nonumber\\[1ex]
&=&\pro\nabla^j\bH\cdot\nabla_j\pro\delta\bH-\nabla^j\bH\cdot\Big[(\bH\cdot\nabla_i\delta\bp)\bh^i_{j}-(\bH\cdot\bh^i_j)\pro\nabla_i\delta\bp  \Big]\nonumber\\[1ex]
&=&\pro\nabla^j\bH\cdot\nabla_j\pro\delta\bH-\nabla^j\bH\cdot\Big[\big[\bH\cdot\nabla_i(A^k\nabla_k\bp+\bB)\big]\bh^i_{j}-(\bH\cdot\bh^i_j)\pro\nabla_i\big(A^k\nabla_k\bp+\bB)  \Big]\nonumber\\[1ex]
&=&\pro\nabla^j\bH\cdot\nabla_j\pro\delta\bH-\nabla^j\bH\cdot\Big[\big(A^k\bH\cdot\bh_{ik}+\bH\cdot\nabla_i\bB\big)\bh^i_j-(\bH\cdot\bh^i_j)\big(A^k\bh_{ik}+\pro\nabla_i\bB\big)  \Big]\nonumber\\[1ex]
&=&\pro\nabla^j\bH\cdot\nabla_j\pro\delta\bH-A^k\Big[(\bH\cdot\bh_{ik})(\bh^i_j\cdot\nabla^j\bH)-(\bH\cdot\bh^i_j)(\bh_{ik}\cdot\nabla^j\bH)\Big]\nonumber\\
&&\hspace{1cm}-\:\big(\bH\cdot\nabla_i\bB)(\bh^i_j\cdot\nabla^j\bH)+(\bH\cdot\bh^i_j)(\nabla^j\bH\cdot\pro\nabla_i\bB)\:.
\end{eqnarray}
On the other hand, one checks that
\begin{eqnarray*}
&&\pro\nabla^j\bH\cdot\nabla_j\pro\nabla_i\bH\;\;=\;\;\pro\nabla^j\bH\cdot\nabla_j\nabla_i\bH-\pro\nabla^j\bH\cdot\nabla_j\pi_T\nabla_i\bH\\[1ex]
&=&\pro\nabla^j\bH\cdot\nabla_i\nabla_j\bH+\pro\nabla^j\bH\cdot\nabla_j\big((\bH\cdot\bh_{is})\nabla^s\bp\big)\\[1ex]
&=&\pro\nabla^j\bH\cdot\nabla_i\pro\nabla_j\bH+\pro\nabla^j\bH\cdot\nabla_i\pi_T\nabla_j\bH+(\bH\cdot\bh_{is})(\bh_j^{s}\cdot\nabla^j\bH)\\[1ex]
&=&\dfrac{1}{2}\nabla_i|\pro\nabla\bH|^2-\pro\nabla^j\bH\cdot\nabla_i\big((\bH\cdot\bh_{sj})\nabla^s\bp\big)+(\bH\cdot\bh_{is})(\bh_j^{s}\cdot\nabla^j\bH)\\[1ex]
&=&\dfrac{1}{2}\nabla_i|\pro\nabla\bH|^2-(\bH\cdot\bh_{sj})(\bh^{s}_i\cdot\nabla^j\bH)+(\bH\cdot\bh_{is})(\bh_j^{s}\cdot\nabla^j\bH)\\[1ex]
&=&\dfrac{1}{2}\nabla_i|\pro\nabla\bH|^2-(\bH\cdot\bh^s_{j})(\bh_{is}\cdot\nabla^j\bH)+(\bH\cdot\bh_{is})(\bh_j^{s}\cdot\nabla^j\bH)\\[1ex]
\end{eqnarray*}
Bringing this in (\ref{eq10}) yields
\begin{eqnarray}\label{eq11}
&&\pro\nabla^j\bH\cdot\delta\pro\nabla_j\bH\;\;=\;\;\pro\nabla^j\bH\cdot\Big[\nabla_j\pro\delta\bH-A^i\nabla_j\pro\nabla_i\bH\Big] +\dfrac{A^i}{2}\nabla_i|\pro\nabla\bH|^2 \nonumber\\
&&\hspace{4cm}-\:\big(\bH\cdot\nabla_i\bB)(\bh^i_j\cdot\nabla^j\bH)+(\bH\cdot\bh^i_j)(\nabla^j\bH\cdot\pro\nabla_i\bB)\:.
\end{eqnarray}
Calling upon (\ref{pindelh}) now gives
\begin{eqnarray*}
&&\pro\nabla^j\bH\cdot\delta\pro\nabla_j\bH\;\;=\;\;\pro\nabla^j\bH\cdot\Big[\dfrac{1}{4}\nabla_j\big(4A^s\pro\nabla_s\bH+(\bB\cdot\bh^{ik})\bh_{ik}+\Delta_\perp\bB\big)-A^i\nabla_j\pro\nabla_i\bH\Big]  \nonumber\\
&&\hspace{4cm}-\:\big(\bH\cdot\nabla_i\bB)(\bh^i_j\cdot\nabla^j\bH)+(\bH\cdot\bh^i_j)(\nabla^j\bH\cdot\pro\nabla_i\bB)+\dfrac{A^i}{2}\nabla_i|\pro\nabla\bH|^2\nonumber\\[1ex]
&&=\:\big(\pro\nabla^j\bH\cdot\pro\nabla_i\bH)\nabla_jA^i+\dfrac{1}{4}\pro\nabla^j\bH\cdot\nabla_j\big((\bB\cdot\bh^{ik})\bh_{ik}\big)+\dfrac{1}{4}\pro\nabla^j\bH\cdot\nabla_j\Delta_\perp\bB\nonumber\\
&&\hspace{4cm}-\:\big(\bH\cdot\nabla_i\bB)(\bh^i_j\cdot\nabla^j\bH)+(\bH\cdot\bh^i_j)(\nabla^j\bH\cdot\pro\nabla_i\bB)+\dfrac{A^i}{2}\nabla_i|\pro\nabla\bH|^2\nonumber\\[1ex]
&&=\:(\bT_i\cdot\bT_j)\nabla^jA^i+\dfrac{A^i}{2}\nabla_i|\bT|^2+\dfrac{1}{4}\bT^j\cdot\nabla_j\big((\bB\cdot\bh^{ik})\bh_{ik}\big)+\dfrac{1}{4}\bT^j\cdot\nabla_j\Delta_\perp\bB\nonumber\\
&&\hspace{4cm}-\:\big(\bH\cdot\nabla_i\bB)(\bh^i_j\cdot\bT^j)+(\bH\cdot\bh^i_j)(\bT^j\cdot\nabla_i\bB)
\end{eqnarray*}
Putting this into (\ref{vary}), it comes
\begin{eqnarray}
&&|g|^{-1/2}\delta\big(|\bT|^2|g|^{1/2}\big)\;\;=\;\;2(\bT_i\cdot\bT_j)(\bB\cdot\bh^{ij})+\dfrac{1}{2}\bT^j\cdot\nabla_j\big((\bB\cdot\bh^{ik})\bh_{ik}\big)+\dfrac{1}{2}\bT^j\cdot\nabla_j\Delta_\perp\bB\nonumber\\
&&\hspace{2cm}-\:2\big(\bH\cdot\nabla_j\bB)(\bh^j_i\cdot\bT^i)+2(\bH\cdot\bh^j_i)(\bT^i\cdot\nabla_j\bB)-4|\bT|^2\bH\cdot\bB+\nabla^j\big(|\bT|^2A_j\big)\nonumber\\[1ex]
&&=\bB\cdot\Big[2(\bT_i\cdot \bT_j)\bh^{ij}-4|\bT|^2\bH-\dfrac{1}{2}(\bh_{ik}\cdot\nabla_j\bT^j)\bh^{ik}+2
\nabla_j\big((\bT^i\cdot\bh^j_i)\bH\big)-2\nabla_j\big((\bH\cdot\bh^j_i)\bT^i\big)\nonumber\\
&&\hspace{12cm} -\:\dfrac{1}{2}\Delta_\perp\pro\nabla_i\bT^i   \Big]     \nonumber\\
&&+\:\nabla^j\Big[|\bT|^2A_j+\dfrac{1}{2}(\bB\cdot\bh^{ik})(\bT_j\cdot\bh_{ik}) -2(\bB\cdot\bH)(\bT^i\cdot\bh_{ij})+2(\bB\cdot\bT^i)(\bH\cdot\bh_{ij})\nonumber\\
&&\hspace{5cm}+\:\dfrac{1}{2}\bT_j\cdot\Delta_\perp\bB-\dfrac{1}{2}\nabla_j\bB\cdot\pro\nabla_i\bT^i+\dfrac{1}{2}\bB\cdot\nabla_j\pro\nabla_i\bT^i \Big]\:.
\end{eqnarray}
This expression can be a little bit simplified by noting that
$$
\pro\nabla_i\bT^i\;=\;\Delta_\perp\bH\:.
$$
Then, we have
$$
|g|^{-1/2}\delta\big(|\bT|^2|g|^{1/2}\big)\;\;=\;\;\bB\cdot\vec{\mathcal{W}}+\nabla^jV_j\:,
$$
with
\begin{eqnarray*}
&&\vec{\mathcal{W}}\;\;=\;\;2(\pro\nabla_i\bH\cdot\pro\nabla_j\bH)\bh^{ij}-4|\pro\nabla\bH|^2\bH-\dfrac{1}{2}(\bh_{ik}\cdot\Delta_\perp\bH)\bh^{ik}\nonumber\\
&&\hspace{3cm}+\:2\pro\nabla_j\big(\bh^j_i\cdot\nabla^i\bH)\bH\big)-2\pro\nabla_j\big((\bH\cdot\bh^j_i)\pro\nabla^i\bH\big) -\dfrac{1}{2}\Delta_\perp^2\bH
\end{eqnarray*}
and
\begin{eqnarray*}
&&V_j\;\;=\;\;|\pro\nabla\bH|^2A_j+\dfrac{1}{2}(\bB\cdot\bh^{ik})(\bh_{ik}\cdot\nabla_j\bH) -2(\bB\cdot\bH)(\bh_{ij}\cdot\nabla^i\bH)\nonumber\\[1ex]
&&+\:2(\bB\cdot\nabla^i\bH)(\bH\cdot\bh_{ij})+\dfrac{1}{2}\pro\nabla_j\bH\cdot\Delta_\perp\bB-\dfrac{1}{2}\nabla_j\bB\cdot\Delta_\perp\bH+\dfrac{1}{2}\bB\cdot\nabla_j\Delta_\perp\bH\:.
\end{eqnarray*}

Suppose next that $\delta\bp=\ba$ is a constant vector in $\mathbb{R}^m$. Then
\be\label{varoo}
A_j\;=\;\ba\cdot\nabla_j\bp\qquad\text{and}\qquad \bB\;=\;\pro\ba\:.
\ee
It is easy to check that if $\bU$ is any vector, it holds
$$
\bU\cdot\bB\;=\;\ba\cdot\pro\bU\qquad\text{and}\qquad\pro\bU\cdot\nabla_j\bB\;=\;\ba\cdot\pi_T\nabla_j\pro\bU\:.
$$
Hence
\begin{eqnarray*}
\bU\cdot\Delta_\perp\bB&=&\pro\bU\cdot\nabla^i\pro\nabla_i\bB\;\;=\;\;\nabla^i\big(\pro\bU\cdot\nabla_i\bB\big)-(\pro\nabla^i\pro\bU)\cdot\nabla_i\bB\\[1ex]
&=&\nabla^i\big[\ba\cdot\pi_T\nabla_i\pro\bU\big]-\ba\cdot\pi_T\nabla_i\pro\nabla^i\pro\bU\\[1ex]
&=&\ba\cdot\Big[\nabla^i\pi_T\nabla_i\pro\bU-\pi_T\nabla_i\pro\nabla^i\pro\bU   \Big]\:.
\end{eqnarray*}
Then we have
$$
V_j\;=\;\ba\cdot\bV_j\:,
$$
with
\begin{eqnarray*}
\bV_j&=&|\pro\nabla\bH|^2\nabla_j\bp+\dfrac{1}{2}(\bh_{ik}\cdot\nabla_j\bH)\bh^{ik} -2(\bh_{ij}\cdot\nabla^i\bH)\bH+2(\bH\cdot\bh_{ij})\pro\nabla^j\bH\nonumber\\
&&+\;\dfrac{1}{2}\pro\nabla_j\Delta_\perp\bH-\dfrac{1}{2}\pi_T\nabla_j\Delta_\perp\bH+\dfrac{1}{2}\nabla^i\pi_T\nabla_i\pro\nabla_j\bH-\dfrac{1}{2}\pi_T\nabla_i\pro\nabla^i\pro\nabla_j\bH\nonumber\\[1ex]
&=&|\pro\nabla\bH|^2\nabla_j\bp-2(\pro\nabla_j\bH\cdot\pro\nabla_i\bH)\nabla^i\bp+\dfrac{1}{2}(\bh_{ij}\cdot\Delta_\perp\bH)\nabla^i\bp+\dfrac{1}{2}\pro\nabla_j\Delta_\perp\bH\nonumber\\
&&-\;2(\bh_{ij}\cdot\nabla^i\bH)\bH+2(\bH\cdot\bh_{ij})\pro\nabla^i\bH\:.
\end{eqnarray*}

\noindent
Variations of $\int|d\bH|^2$ yield the equation $d^\star\bV=-\vec{\mathcal{W}}$. At equilibrium, for critical points, the term $\vec{\mathcal{W}}$ vanishes.  \\

We will find convenient to write
$$
\bV^j\;=\;G^{ij}\nabla_i\bp-\pro\nabla_i\bF^{ij}+\bC^j\:,
$$
where
\be\label{GAFA}
\left\{\begin{array}{lcl}
G^{ij}&=&\dfrac{1}{2}\bh^{ij}\cdot\Delta_\perp\bH+ |\pro\nabla\bH|^2g^{ij}-2\nabla_j\bH\cdot\pro\nabla_i\bH\\[1ex]
\bF^{ij}&=&-\dfrac{1}{2}\Delta_\perp\bH\,g^{ij} \\[1ex]
\bC^j&=&-2(\bh^{ij}\cdot\nabla_i\bH)\bH+2(\bH\cdot\bh^{ij})\pro\nabla_i\bH\:.\end{array}\right.
\ee

\begin{Rm}\label{rem0}
The tensor $G$ satisfies
$$
G^i_i\;=\;2\bH\cdot\Delta_\perp\bH+2|\pro d\bH|^2\;=\;\Delta|\bH|^2\:.
$$
\end{Rm}
$\hfill\square$\\

The triple of data $(G,\bF,\bC)$ are of course closely tied to each other. For our use, we record in particular the following fact. 
\begin{Lm}\label{lem0}
It holds
$$
\bh^m_j\cdot\nabla_i\bF^{ij}-\bh^m_j\cdot\bC^j\;=\;-\nabla_iG^{im}\:.
$$
\end{Lm}
{\bf Proof.}
With the help of Codazzi, it holds
\begin{eqnarray*}
-\bh^m_j\cdot\nabla_i\bF^{ij}+\bh^m_j\cdot\bC^j&=&-\nabla_i\big(\bh^m_j\cdot\bF^{ij}\big)+\bF^{ij}\cdot\nabla^m\bh_{ij}+\bh^m_j\cdot\bC^j\\
&=&-\nabla_i\big(\bh^m_j\cdot\bF^{ij}-g^{mi}\langle\bh\stackrel{\cdot}{,}\bF\rangle \big)-\bh_{ij}\cdot\nabla^m\bF^{ij}+\bh^m_j\cdot\bC^j\:.
\end{eqnarray*}
Since 
$$
\bh^m_j\cdot\bF^{ij}-g^{mi}\langle\bh\stackrel{\cdot}{,}\bF\rangle\;=\;-G^{im}+g^{im}\big|\pro d\bH\big|^2-2\nabla^i\bH\cdot\pro\nabla^m\bH+2\bH\cdot\Delta_\perp\bH\:.
$$
the latter yields
\begin{eqnarray}\label{PS1}
-\bh^m_j\cdot\nabla_i\bF^{ij}+\bh^m_j\cdot\bC^j&=&-\bh_{ij}\cdot\nabla^m\bF^{ij}+\bh^m_j\cdot\bC^j\nonumber\\
&&\hspace{-1.5cm}+\:\nabla_i\big(G^{im}-g^{im}\big|\pro d\bH\big|^2+2\nabla^i\bH\cdot\pro\nabla^m\bH-2\bH\cdot\Delta_\perp\bH\big)\:.
\end{eqnarray}
On the other hand, we have
\begin{eqnarray*}
-\bh_{ij}\cdot\nabla^m\bF^{ij}&=&2\bH\cdot\nabla^m\Delta_\perp\bH\\
&=&2\nabla^m\big(\bH\cdot\Delta_\perp\bH\big)-2\Delta_\perp\bH\cdot\nabla^m\bH\\
&=&2\nabla^m\big(\bH\cdot\Delta_\perp\bH\big)-2\nabla_j\big(\pro\nabla^j\bH\cdot\pro\nabla^m\bH\big)\\
&&\hspace{1cm}+\:2\pro\nabla_j\bH\cdot\nabla^j\pro\nabla^m\bH\\
&=&2\nabla^m\big(\bH\cdot\Delta_\perp\bH\big)-2\nabla_j\big(\pro\nabla^j\bH\cdot\pro\nabla^m\bH\big)+\nabla^m\big|\pro d\bH\big|^2\\
&&\hspace{1cm}+\:2\pro\nabla_j\bH\cdot\big(\nabla^m\pi_T\nabla^j\bH-\nabla^j\pi_T\nabla^m\bH\big)\\
&=&2\nabla^m\big(\bH\cdot\Delta_\perp\bH\big)-2\nabla_j\big(\pro\nabla^j\bH\cdot\pro\nabla^m\bH\big)+\nabla^m\big|\pro d\bH\big|^2\\
&&\hspace{1cm}+\:2\big((\bH\cdot\bh^{m}_{j})\bh^{ij}-(\bH\cdot\bh^{ij})\bh^m_j\big)\cdot\nabla_i\bH\:.
\end{eqnarray*}
Equivalently, the latter states that
$$
-\bh_{ij}\cdot\nabla^m\bF^{ij}+\bh^m_j\cdot\bC^j\;=\;\nabla_j\big(2g^{mj}\bH\cdot\Delta_\perp\bH+g^{mj}|\pro d\bH|^2-2\nabla^j\bH\cdot\pro\nabla^m\bH \big)\:.
$$
Importing this into (\ref{PS1}) confirms that
\be\label{PS2}
-\bh^m_j\cdot\nabla_i\bF^{ij}+\bh^m_j\cdot\bC^j\;=\;\nabla_iG^{im}\:.
\ee
$\hfill\blacksquare$

The next result will also be useful in our derivations. 

\begin{Lm}\label{lem4}
It holds
$$
\big(-\pro\nabla_i\bF^{ij}+\bC^j\big)\wedge\nabla_j\bp\;=\;\nabla_i\big(-\bF^{ij}\wedge\nabla_j\bp+2\bH\wedge\pro\nabla^i\bH\big)\:.
$$
\end{Lm}
{\bf Proof.} 
Define the one-vector-valued one-form
$$
\bX^j\;:=\;\pro\nabla^j\bH-\big(2|\bH|^2g^{jk}-\bH\cdot\bh^{jk}\big)\nabla_k\bp\:.
$$
We verify the exact divergence
\begin{eqnarray*}
\nabla_j\bX^j&=&\Delta_\perp\bH-(\bh^{jk}\cdot\nabla_j\bH)\nabla_k\bp-2\nabla_j|\bH|^2\nabla^j\bp+(\bh^{jk}\cdot\nabla_j\bH)\nabla_k\bp+2\nabla_j|\bH|^2\nabla^j\bp\\
&&\hspace{1cm}-\:\big(2|\bH|^2g^{jk}-\bH\cdot\bh^{jk}\big)\bh_{jk}\\
&=&\Delta_\perp\bH-\big(2|\bH|^2g^{jk}-\bH\cdot\bh^{jk}\big)\bh_{jk}\:.
\end{eqnarray*}
Accordingly, we have
\begin{eqnarray*}
\bH\wedge\big(\Delta_\perp\bH+(\bH\cdot\bh^{ij})\bh_{ij}\big)&=&\bH\wedge\nabla_j\bX^j\\
&=&\nabla_j(\bH\wedge\bX^j)+\bX^j\wedge\nabla_j\bH\\
&=&\nabla_j(\bH\wedge\bX^j)+\pro\nabla^j\bH\wedge\pi_T\nabla_j\bH+2|\bH|^2\nabla^j\bH\wedge\nabla_j\bp\\
&&\hspace{1cm}-\:(\bh^{ij}\cdot\bH)\nabla_i\bH\wedge\nabla_j\bp\\
&=&\nabla_j(\bH\wedge\bX^j)+2|\bH|^2\nabla^j\bH\wedge\nabla_j\bp-2(\bh^{ij}\cdot\bH)\nabla_i\bH\wedge\nabla_j\bp\:.
\end{eqnarray*}
Next, using Codazzi, we find
\begin{eqnarray*}
&&\nabla_i\big((\bh^{ij}\cdot\bH)\bH\wedge\nabla_j\bp\big)\nonumber\\[1ex]
&&\hspace{1cm}=\:(\bh^{ij}\cdot\bH)\nabla_i\bH\wedge\nabla_j\bp+(\bh^{ij}\cdot\bH)\bH\wedge\bh_{ij}+2\nabla^j|\bH|^2\bH\wedge\nabla_j\bp\nonumber\\
&&\hspace{2cm}+\:(\bh^{ij}\cdot\nabla_i\bH)\bH\wedge\nabla_j\bp\:.
\end{eqnarray*}
Combining together the last two identities and reorganising the terms yields
\begin{eqnarray}\label{PS5}
\bH\wedge\Delta_\perp\bH&=&\nabla_j(\bH\wedge\bX^j)+2|\bH|^2\nabla^j\bH\wedge\nabla_j\bp-2(\bh^{ij}\cdot\bH)\nabla_i\bH\wedge\nabla_j\bp-(\bH\cdot\bh^{ij})\bH\wedge\nabla\bh_{ij}\nonumber\\
&=&\nabla_j(\bH\wedge\bX^j)+2|\bH|^2\nabla^j\bH\wedge\nabla_j\bp-(\bh^{ij}\cdot\bH)\nabla_i\bH\wedge\nabla_j\bp+2\nabla_j|\bH|^2\bH\wedge\nabla_j\bp\nonumber\\
&&\hspace{1cm}+\:\big(\bh^{ij}\cdot\nabla_i\bH)\bH\wedge\nabla_j\bp-\nabla_i\big((\bh^{ij}\cdot\bH)\bH\wedge\nabla_j\bp\big)\nonumber\\
&=&\nabla_j\big(\bH\wedge\bX^j+2|\bH|^2\bH\wedge\nabla^j\bp-(\bh^{ij}\cdot\bH)\bH\wedge\nabla_i\bp\big)\nonumber\\
&&\hspace{1cm}+\:\big((\bh^{ij}\cdot\nabla_i\bH)\bH-(\bh^{ij}\cdot\bH)\pro\nabla_i\bH\big)\wedge\nabla_j\bp\:.
\end{eqnarray}
Note that we have used the fact that
$$
(\bh^{ij}\cdot\bH)\pi_T\nabla_i\bH\wedge\nabla_j\bp\;=\;-(\bh^{ij}\cdot\bH)(\bh_{i}^k\cdot\bH)\nabla_k\bp\wedge\nabla_j\bp\;=\;\vec{0}\:.
$$
Without much effort now, we call upon (\ref{GAFA}) to recast (\ref{PS5}) in the equivalent form
\begin{eqnarray*}\label{PS6}
\bF^{ij}\wedge\bh_{ij}&=&-2\Delta_\perp\bH\wedge\bH\;\;=\;\;2\bH\wedge\Delta_\perp\bH\\
&=&2\nabla_j\big(\bH\wedge\bX^j+2|\bH|^2\bH\wedge\nabla^j\bp-(\bh^{ij}\cdot\bH)\bH\wedge\nabla_i\bp\big)-\bC^j\wedge\nabla_j\bp\:.
\end{eqnarray*}
Finally, since $\bF^{ij}\cdot\bh_i^k$ is clearly a symmetric 2-tensor, it comes
\begin{eqnarray*}
\pro\nabla_i\bF^{ij}\wedge\nabla_j\bp&=&\nabla_i(\bF^{ij}\wedge\nabla_j\bp)-\bF^{ij}\wedge\bh_{ij}+(\bF^{ij}\cdot\bh_{i}^k)\nabla_k\bp\wedge\nabla_j\bp\nonumber\\
&=&\nabla_i(\bF^{ij}\wedge\nabla_j\bp)-\bF^{ij}\wedge\bh_{ij}\:,
\end{eqnarray*}
which, once combined with (\ref{PS6}), implies the announced
\begin{eqnarray*}
&&-\pro\nabla_i\bF^{ij}\wedge\nabla_j\bp+\bC^j\wedge\nabla_j\bp\nonumber\\
&&\hspace{1cm}\:=\:\nabla_i\big(-\bF^{ij}\wedge\nabla_j\bp+2\bH\wedge\bX^i+4|\bH|^2\bH\wedge\nabla^i\bp-2(\bh^{ij}\cdot\bH)\bH\wedge\nabla_i\bp\big)\\
&&\hspace{1cm}\:=\:\nabla_i\big(-\bF^{ij}\wedge\nabla_j\bp+2\bH\wedge\pro\nabla^i\bH\big)\:.
\end{eqnarray*}
$\hfill\blacksquare$

\subsubsection{Variation of $\int\mathcal{O}(|\bh|^4)$}

In dimension 4, there are exactly 8 distinct energies of the type $\int\mathcal{O}(|\bh|^4)$. Letting $T_{ab}^{cd}:=\bh_a^c\cdot\bh_b^d$, they can all eight be obtained from various contractions of the quantity $T_{ab}^{cd}T^{kl}_{ij}$. That is therefore the quantity which we will vary. We proceed exactly as in the previous section and consider a variation of the type (\ref{varo}). We open with
\begin{eqnarray*}
\pro\delta\bh_a^c&=&\pro\delta(g^{ci}\bh_{ai})\\[1ex]
&=&g^{ci}\pro\delta\bh_{ai}+\bh_{ai}\delta g^{ci}\\[1ex]
&=&g^{ci}\big(A^s\pro\nabla_{a}\bh_{si}+\bh_{sa}\nabla_iA^s+\bh_{si}\nabla_aA^s+\pro\nabla_{ai}\bB\big)-\bh_{ai}\big(\nabla^c A^i+\nabla^i A^c-2\bB\cdot\bh^{ic}  \big)\:.
\end{eqnarray*}
From this, it follows that
\begin{eqnarray*}
\delta T_{ab}^{cd}&=&\bh_a^c\cdot\pro\delta\bh^d_b+\bh_b^d\cdot\pro\delta\bh^c_a\\[1ex]
&=&\bh^{d}_{b}\cdot\Big[\big(A^s\pro\nabla_{a}\bh_{s}^c+\bh^s_{a}\nabla^cA_s+\bh_s^{c}\nabla_aA^s+\pro\nabla^c_{a}\bB\big)-\bh_{a}^{i}\big(\nabla^c A
_i+\nabla_i A^c-2\bB\cdot\bh_{i}^{c}  \big)\Big]\\
&&+\:\bh^{c}_{a}\cdot\Big[\big(A^s\pro\nabla_{b}\bh_{s}^d+\bh^s_{b}\nabla^dA_s+\bh^{d}_s\nabla_bA^s+\pro\nabla^d_{b}\bB\big)-\bh_{b}^{i}\big(\nabla^d A
_i+\nabla_i A^d-2\bB\cdot\bh_{i}^{d}  \big)\Big]\\[1ex]
&=&A^s\bh^d_b\cdot\nabla_s\bh^c_a+T_{ba}^{ds}\nabla^cA_s+T_{bs}^{dc}\nabla_aA^s+\bh_b^d\cdot\nabla_a^c\bB-T_{ba}^{ds}\big(\nabla^cA_s+\nabla_sA^c-2\bB\cdot\bh^{c}_{s}\big)\\
&&+\:A^s\bh^c_a\cdot\nabla_s\bh^d_b+T_{ab}^{cs}\nabla^dA_s+T_{as}^{cd}\nabla_bA^s+\bh_a^c\cdot\nabla_b^d\bB-T_{ab}^{cs}\big(\nabla^dA_s+\nabla_sA^d-2\bB\cdot\bh^{d}_{s}\big)\\[1ex]
&=&A^s\nabla_sT^{cd}_{ab}+2(\bB\cdot\bh^c_s)T^{ds}_{ba}+2(\bB\cdot\bh^d_s)T^{cs}_{ab}\\
&&+\:T_{bs}^{dc}\nabla_aA^s+T_{as}^{cd}\nabla_bA^s-T^{ds}_{ba}\nabla_sA^c-T^{cs}_{ab}\nabla_sA^d+\bh_b^d\cdot\nabla_a^c\bB+\bh^c_a\cdot\nabla_b^d\bB\\[1ex]
&=&A^s\nabla_sT^{cd}_{ab}+2(\bB\cdot\bh^c_s)T^{ds}_{ba}+2(\bB\cdot\bh^d_s)T^{cs}_{ab}+\bB\cdot\big(\ldots  \big)\\
&&+\:T_{bs}^{dc}\nabla_aA^s+T_{as}^{cd}\nabla_bA^s-T^{ds}_{ba}\nabla_sA^c-T^{cs}_{ab}\nabla_sA^d+\bh_b^d\cdot\nabla_a^c\bB+\bh^c_a\cdot\nabla_b^d\bB\\[1ex]
&=&A^s\nabla_sT^{cd}_{ab}+\bB\cdot(\bh^c_sT^{ds}_{ba}+\bh_a^sT^{dc}_{bs}+\bh_b^sT^{cd}_{as}+\bh^d_sT^{cs}_{ab})+\bh_b^d\cdot\nabla_a^c\bB+\bh^c_a\cdot\nabla_b^d\bB+\bB\cdot\big(\ldots  \big)\:.
\end{eqnarray*}
We have deliberately ignored any term of the type $\bB\cdot(\ldots)$ as it gives rise to an identically vanishing term at equilibrium. Our computations imply
\begin{eqnarray*}
\delta(T^{cd}_{ab}T_{ij}^{kl})&=&A^s\nabla_s(T^{cd}_{ab}T_{ij}^{kl})+(\bh^d_bT_{ij}^{kl})\cdot\nabla^c_a\bB+(\bh^c_aT_{ij}^{kl})\cdot\nabla^d_b\bB\\
&&\hspace{1cm}+\:(\bh^k_iT_{ab}^{cd})\cdot\nabla^l_j\bB+(\bh^l_jT_{ab}^{cd})\cdot\nabla^k_i\bB+\bB\cdot\big(\ldots  \big)\:,
\end{eqnarray*}
and thus
\begin{eqnarray*}
|g|^{-1/2}\delta\big(|g|^{1/2}T^{cd}_{ab}T_{ij}^{kl}\big)&=&\nabla^m\big(T_{ab}^{cd}T_{ij}^{kl}A_m\big)+\bB\cdot\big(\ldots  \big)\\
&&\hspace{-1cm}+\;(\bh^d_bT_{ij}^{kl})\cdot\nabla^c_a\bB+(\bh^c_aT_{ij}^{kl})\cdot\nabla^d_b\bB+(\bh^k_iT_{ab}^{cd})\cdot\nabla^l_j\bB+(\bh^l_jT_{ab}^{cd})\cdot\nabla^k_i\bB\:.
\end{eqnarray*}
Identically as in the previous section, we now specialise to a constant translation $\delta\bp=\vec{a}$, as in (\ref{varoo}). Noting that for any normal vector $\bU$, it holds
\begin{eqnarray*}
\bU\cdot\nabla^p_q\bB&\simeq&\vec{a}\cdot\Big(\pro\nabla_q^{\;\;p}\bU+\nabla^m\big[\delta^p_s\pi_{\vec{t}}\nabla_q\bU-g_{qs}\pro\nabla^p\bU  \big]\Big)\\[1ex]
&\simeq&-\vec{a}\cdot\nabla^s\big((\delta^p_s\bh_{qr}\cdot\bU\big)\nabla^r\bp+g_{qs}\pro\nabla^p\bU  \big)\:,
\end{eqnarray*}
we can now derive
\begin{eqnarray*}
&&|g|^{-1/2}\delta\big(|g|^{1/2}T^{cd}_{ab}T_{ij}^{kl}\big)\;\;=\;\;\vec{a}\cdot\pro\big(\ldots\big)\\[1ex]
&&\hspace{-1cm}+\;\vec{a}\cdot\nabla^s\Big[\big(T^{cd}_{ab}T_{ij}^{kl}g_{rs}
-\delta^c_s\bh_{ar}\cdot\bh^d_bT^{kl}_{ij}-\delta^d_s\bh_{br}\cdot h^{c}_{a}T_{ij}^{kl}-\delta^l_s\bh_{jr}\cdot\bh^k_iT_{ab}^{cd}-\delta^k_s\bh_{ir}\cdot\bh^l_jT_{ab}^{cd}\big)\nabla^r\bp\\
&&\hspace{1.5cm}-g_{as}\pro\nabla^c\big(\bh^d_bT^{kl}_{ij}\big)-g_{bs}\pro\nabla^d\big(\bh_a^cT^{kl}_{ij}\big)-g_{js}\pro\nabla^l\big(\bh^k_iT^{cd}_{ab}\big)-g_{is}\pro\nabla^k\big(\bh^l_jT^{cd}_{ab}\big)\Big]\\[1ex]
&&\hspace{-1cm}=\;\vec{a}\cdot\nabla^s\Big[\big(\delta^{r}_{s}T^{cd}_{ab}T_{ij}^{kl}
-\delta^c_s T^{rd}_{ab}T^{kl}_{ij}-\delta^d_sT^{cr}_{ab}T_{ij}^{kl}-\delta^l_s T^{kr}_{ij}T_{ab}^{cd}-\delta^k_s T^{rl}_{ij}T_{ab}^{cd}\big)\nabla_r\bp\\
&&\hspace{1cm}-\pro\nabla_r\big(g_{as}g^{rc}\bh^d_bT^{kl}_{ij}+g_{bs}g^{dr}\bh^c_aT^{kl}_{ij}+g_{js}g^{lr}\bh^k_iT^{cd}_{ab}+g_{is}g^{kr}\bh^l_jT^{cd}_{ab}\big)\Big]+\vec{a}\cdot\pro\big(\ldots\big)\:.\\[1ex]
\end{eqnarray*}
We let
\begin{equation}\label{defGF}
\left\{\begin{array}{lcl}
G^{r}_s&:=&\delta^{r}_{s}T^{cd}_{ab}T_{ij}^{kl}
-\delta^c_s T^{rd}_{ab}T^{kl}_{ij}-\delta^d_sT^{cr}_{ab}T_{ij}^{kl}-\delta^l_s T^{kr}_{ij}T_{ab}^{cd}-\delta^k_s T^{rl}_{ij}T_{ab}^{cd}\\[1ex]
\bF^r_s&:=&g_{as}g^{rc}\bh^d_bT^{kl}_{ij}+g_{bs}g^{dr}\bh^c_aT^{kl}_{ij}+g_{js}g^{lr}\bh^k_iT^{cd}_{ab}+g_{is}g^{kr}\bh^l_jT^{cd}_{ab}\:,
\end{array}\right.
\end{equation}
and write 
$$
|g|^{-1/2}\delta\big(|g|^{1/2}T^{cd}_{ab}T_{ij}^{kl}\big)\;=\;\vec{a}\cdot\Big[\vec{\mathcal{W}}+\nabla^s\big(G^r_s\nabla_r\bp-\pro\nabla_r\bF^r_s\big)\Big]\:,
$$
where $\vec{\mathcal{W}}$ entails all terms previously rendered as $(\ldots)$. Critical points thus satisfy $d^\star\bV=-\vec{\mathcal{W}}$, where
$$
\bV^s\;:=\;G^{rs}\nabla_r\bp-\pro\nabla_r\bF^{rs}\:.
$$

\bigskip

With the help of the above algorithm, we compute the characteristic tensor $\bF$ associated for each of the eight fundamental energies with integrand of order $\mathcal{O}(|\bh|^4)$. The results are compiled in the following table. 
\begin{center}
\begin{equation}\label{tablebraid}
\begin{tabular}{|c|c|}
\hline
 &\\[-1.2ex]
Energy $\int E$  & $\bF^{rs}$\\[.8ex]
\hline
  &\\[-1.2ex]
 $T^{cd}_{ab}T^{ab}_{cd}\;=\;(\bh^{ij}\cdot\bh^{kl})\bh_{ij}\cdot\bh_{kl}$ & $4(\bh^{rs}\cdot\bh^{ab})\bh_{ab}$\\[1ex]
\hline
  &\\[-1.2ex]
 $T^{cd}_{ab}T^{ab}_{dc}=:\text{Tr}_g(\bh^4)=(\bh^{ij}\cdot\bh^{kl})\bh_{ik}\cdot\bh_{jl}$ & $4(\bh^{ra}\cdot\bh^{sb})\bh_{ab}$\\[1ex]
\hline
  &\\[-1.2ex]
 $T^{ca}_{ab}T^{bd}_{dc}\;=\;|h^2|^2$ & $2(h^2)^{rb}\bh^{s}_{b}+2(h^2)^{sb}\bh^r_b$\\[1ex]
\hline
  &\\[-1.2ex]
 $T^{ab}_{ba}T_{cd}^{dc}\;=\;|\bh|^4$ & $4|\bh|^2\bh^{rs}$\\[1ex]
\hline
  &\\[-1.2ex]
 $T^{ab}_{ac}T^{ci}_{ib}\;=\;4\langle h^2,\bH\cdot\bh\rangle$ & $4(\bH\cdot\bh^{rb})\bh^{s}_b+4(\bH\cdot\bh^{sb})\bh^{r}_b+4(h^2)^{rs}\bH+\langle h^2,\bh\rangle g^{rs}$\\[1.75ex]
\hline
 &\\[-1.2ex]
 $T^{ib}_{ic}T^{ac}_{ab}\;=\;4^2|\bH\cdot\bh|^2$ & $32(\bH\cdot\bh^{rs})\bH+8\langle \bH\cdot\bh,\bh\rangle g^{rs}$\\[1.75ex]
\hline
  &\\[-1.2ex]
 $T^{ab}_{ab}T^{cd}_{dc}\;=\;4^2|\bH|^2|\bh|^2$ & $32|\bH|^2\bh^{rs}+8|\bh|^2\bH g^{rs}$\\[1.75ex]
\hline
  &\\[-1.2ex]
 $T^{ab}_{ab}T^{cd}_{cd}\;=\;4^4|\bH|^4$ & $4^4|\bH|^2\bH g^{r}_s$\\[-1ex]
   &\\
\hline
\end{tabular}
\end{equation}
\end{center} 
For notational convenience, we have taken the liberty to define the symmetric 2-tensor
$$
(h^2)^{ij}\;:=\;\bh^{ik}\cdot\bh_k^j\:.
$$
\begin{Rm}\label{rem1}
We have only listed the $\bF$ tensors, as the $G$ tensors can be obtained from $\bF$ as follows:
$$
G^{rs}\;=\;-\bF^{rb}\cdot\bh^{s}_b+Eg^{rs}\:.
$$
This is an easy consequence of (\ref{defGF}).\\
From the table, we note that 
$$
\bh^{ik}\cdot\bF_k^j\;=\;\bh^{jk}\cdot\bF_k^i\qquad\text{and}\qquad\langle\bh\stackrel{\cdot}{,}\bF\rangle\;=\;4E\:.
$$
Hence $G$ must be traceless. 
\end{Rm}
$\hfill\square$\\

The data given in the table will be used in the next section. We also need to record the following observation. 
\begin{Lm}\label{lem3}
For an energy with integrand 
$$
E\;=\;|\bH\cdot\bh|^2+\alpha|\bH|^4\:\qquad\text{where $\alpha$ is a constant},
$$
it holds
$$
\bF^{rs}\wedge\bh_{rs}\;=\;\vec{0}\qquad\text{and}\qquad\bh^m_s\cdot\nabla_r\bF^{rs}\;=\;-\nabla_rG^{rs}\:.
$$
\end{Lm}
{\bf Proof.} 
Per the table,
$$
\bF^{rs}\;=\;2(\bH\cdot\bh^{rs})\bH+\dfrac{1}{2}\langle\bH\cdot\bh,\bh\rangle g^{rs}+\al|\bH|^2\bH\,g^{rs}\:.
$$
The first item follows immediately. To reach the second item, we observe that
\begin{eqnarray*}
\nabla^m|\bH\cdot\bh|^2&=&2(\bH\cdot\bh_{rs})\nabla^m(\bH\cdot\bh^{rs})\\
&=&2\bh_{rs}\cdot\nabla^m\big((\bH\cdot\bh^{rs})\bH\big)-2\langle\bH\cdot\bh,\bh\rangle\cdot\nabla^m\bH\\
&=&2\bh_{rs}\cdot\nabla^m\big((\bH\cdot\bh^{rs})\bH\big)-2\nabla^m\big(\langle\bH\cdot\bh,\bh\rangle\cdot\bH \big)+\dfrac{1}{2}\bh_{rs}\cdot\nabla^m\big(\langle\bH\cdot\bh,\bh\rangle g^{rs}\big)\\
&=&\bh_{rs}\cdot\nabla^m\Big(2(\bH\cdot\bh^{rs})\bH+\dfrac{1}{2}\langle\bH\cdot\bh,\bh\rangle g^{rs}\Big)-2\nabla^m|\bH\cdot\bh|^2\:.
\end{eqnarray*}
Similarly, it holds
\begin{eqnarray*}
\nabla^m|\bH|^4&=&2|\bH|^2\nabla^m|\bH|^2\;\;=\;\;2\bH\cdot\nabla^m\big(|\bH|^2\bH\big)-2|\bH|^2\bH\cdot\nabla^m\bH\\
&=&\dfrac{1}{2}\bh_{rs}\cdot\nabla^m\big(|\bH|^2\bH g^{rs}\big)-\dfrac{1}{2}\nabla^m|\bH|^4\:.
\end{eqnarray*}
From this it easily follows that
\be\label{PS4}
\bh_{rs}\cdot\nabla^m\bF^{rs}\;=\;3\nabla^iE\:.
\ee
To complete the proof, it remains to note that Codazzi implies
\begin{eqnarray*}
\bh^m_s\cdot\nabla_r\bF^{rs}&=&\nabla_r\big(\bh^m_s\cdot\bF^{rs}-\langle\bh\stackrel{\cdot}{,}\bF\rangle g^{mr}\big)+\bh_{rs}\cdot\nabla^m\bF^{rs}\\
&\stackrel{\text{(\ref{PS4})}}{=}&\nabla_r\big(\bh^m_s\cdot\bF^{rs}-\langle\bh\stackrel{\cdot}{,}\bF\rangle g^{mr}\big)+3\nabla^mE\\
&=&-\nabla_rG^{rm}\:,
\end{eqnarray*}
where in the last step we've used Remark~\ref{rem1}.

$\hfill\blacksquare$

\subsubsection{Conformally invariant energies}\label{derder}

We will be interested in conformal invariant energies, thus in those involving the traceless second fundamental form $\bh_0:=\bh-g\bH$. 
A few somewhat tedious, but elementary, computations help convert the data from the previous table into the following ones:
\begin{center}
\begin{equation}\label{tablebraid0}
\begin{tabular}{|c|c|}
\hline
 &\\[-1.2ex]
Energy $\int E$  & $\bF^{rs}$\\[.8ex]
\hline
  &\\[-1.2ex]
 $\langle\bh_0^4\rangle:=(\bh_0^{ij}\cdot\bh_0^{kl})(\bh_0)_{ij}\cdot(\bh_0)_{kl}$ & $4(\bh_0^{rs}\cdot\bh_0^{ab})(\bh_0)_{ab}$\\[1ex]
\hline
  &\\[-1.2ex]
 $\text{Tr}_g(\bh_0^4):=(\bh_0^{ij}\cdot\bh_0^{kl})(\bh_0)_{ik}\cdot(\bh_0)_{jl}$ & $4(\bh_0^{ra}\cdot\bh_0^{sb})(\bh_0)_{ab}-\langle h_0^2,\bh_0\rangle g^{rs}$\\[1ex]
\hline
  &\\[-1.2ex]
 $|h_0^2|^2$ & $2(h_0^2)^{rb}(\bh_0)^{s}_{b}+2(h_0^2)^{sb}(\bh_0)^r_b-\langle h_0^2,\bh_0\rangle g^{rs}$\\[1ex]
\hline
  &\\[-1.2ex]
 $|\bh_0|^4$ & $4|\bh_0|^2\bh_0^{rs}$\\[1ex]
\hline
\end{tabular}
\end{equation}
\end{center} 
Of course, in codimension 1, the first and last energies are the same, while the second and third ones also coincide. We record the following observations:

\begin{Lm}\label{lem1}
For each of the four conformally invariant energies $\int E$, it holds:
$$
\bF^r_r\;=\;\vec{0}\qquad\text{and}\qquad\bF^{rs}\wedge\bh_{rs}\;=\;\vec{0}\:.
$$
In addition, for all indices $m$:
$$
(\bh_0)_{rs}\cdot\nabla^m\bF^{rs}\;=\;3\nabla^mE\:.
$$

\end{Lm}
{\bf Proof.} The first item is evident and follow easily from the tracelessness of $\bh_0$. The second item follows from the tracelessness of $\bF$, so that $\bF^{rs}\wedge\bh_{rs}=\bF^{rs}\wedge(\bh_0)_{rs}$. This quantity can be read off of the data in the table: it vanishes. The third item deserves some attention. We begin with the last energy in the table, namely:
\begin{eqnarray*}
(\bh_0)_{rs}\cdot\nabla^m\bF^{rs}&=&4(\bh_0)_{rs}\cdot\big(|\bh_0|^2\nabla^m\bh_0^{rs}+\bh_0^{rs}\nabla^m|\bh_0|^2\big)\\
&=&2|\bh_0|^2\nabla^m|\bh_0|^2+4|\bh_0|^2\nabla^m|\bh_0|^2\\
&=&3\nabla^m|\bh_0|^4\:,
\end{eqnarray*}
as claimed. \\
Next, for the first energy
\begin{eqnarray*}
(\bh_0)_{rs}\cdot\nabla^m\bF^{rs}&=&\nabla^m\langle\bh_0\stackrel{\cdot}{,}\bF\rangle-\bF^{rs}\cdot\nabla^m(\bh_0)_{rs}\\
&=&4\nabla^m|\bh_0^{ij}\cdot\bh_0^{kl}|^2-4(\bh_0^{rs}\cdot\bh_0^{ab})(\bh_0)_{ab}\cdot\nabla^m(\bh_0)_{rs}\\
&=&4\nabla^m|\bh_0^{ij}\cdot\bh_0^{kl}|^2-2(\bh_0^{rs}\cdot\bh_0^{ab})\big((\bh_0)_{ab}\cdot\nabla^m(\bh_0)_{rs}+(\bh_0)_{rs}\cdot\nabla^m(\bh_0)_{ab}\big)\\
&=&4\nabla^m|\bh_0^{ij}\cdot\bh_0^{kl}|^2-2(\bh_0^{rs}\cdot\bh_0^{ab})\nabla^m\big((\bh_0)_{ab}\cdot(\bh_0)_{rs}\big)\\
&=&3\nabla^m|\bh_0^{ij}\cdot\bh_0^{kl}|^2\:.
\end{eqnarray*}
The third energy gives
\begin{eqnarray*}
2\nabla^m|h_0^2|^2&=&\nabla^m\big(\bh_0^{ik}\cdot(h_0^2)_{ij}(\bh_0)^j_k\big)+\nabla^m\big(\bh_0^{ik}\cdot(h_0^2)_{kj}(\bh_0)^j_i\big)\\
&=&\bh_0^{ik}\cdot\nabla^m\big((h_0^2)_{ij}(\bh_0)^j_k+(h_0^2)_{kj}(\bh_0)^j_i\big)+2(h_0^2)_{ij}(\bh_0)^j_k\cdot\nabla^m\bh_0^{ik}\\
&=&\dfrac{1}{2}(\bh_0)_{ik}\cdot\nabla^m\bF^{ik}+(h_0^2)_{ij}\nabla^m(h_0^2)^{ij}\\
&=&\dfrac{1}{2}(\bh_0)_{ik}\cdot\nabla^m\bF^{ik}+\dfrac{1}{2}\nabla^m|h_0^2|^2\:,
\end{eqnarray*}
which proves the claim. \\
The second energy is handled {\it mutatis mutandis.}

$\hfill\blacksquare$

We deduce from this latter statement two useful results. 

\begin{Co}\label{cor1}
We have
$$
\bh^i_j\cdot\nabla_k\bF^{kj}\;=\;-\nabla_kG^{ik}\:.
$$
\end{Co}
{\bf Proof.} Clearly, using Codazzi,
\begin{eqnarray*}
\bh^i_j\cdot\nabla_k\bF^{kj}&=&\nabla_k\big(\bh^i_j\cdot\bF^{kj}-\langle\bh\stackrel{\cdot}{,}\bF\rangle g^{ik}\big)+\bh_{kj}\cdot\nabla^i\bF^{kj}\:.
\end{eqnarray*}
From Lemma~\ref{lem1}, we know that 
$$
\bh_{kj}\cdot\nabla^i\bF^{kj}\;=\;(\bh_0)_{kj}\cdot\nabla^i\bF^{kj}\;=\;3\nabla^iE\:.
$$
Also, from Remark~\ref{rem1}
$$
\bh^i_j\cdot\bF^{kj}-\langle\bh\stackrel{\cdot}{,}\bF\rangle g^{ik}\;=\;-G^{ik}-3Eg^{ik}\:.
$$
Accordingly now, the desired claim follows:
$$
\bh^i_j\cdot\nabla_k\bF^{kj}\;=\;-\nabla_kG^{ik}\:.
$$
$\hfill\blacksquare$

\begin{Co}\label{cor2}
We have
$$
-\pro\nabla_i\bF^{ij}\wedge\nabla_j\bp\;=\;-\nabla_i\big(\bF^{ij}\wedge\nabla_j\bp\big)\:.
$$
\end{Co}
{\bf Proof.} This is immediate from Lemma~\ref{lem1}. Indeed,
\begin{eqnarray*}
-\pro\nabla_i\bF^{ij}\wedge\nabla_j\bp&=&-\nabla_i\big(\bF^{ij}\wedge\nabla_j\bp\big)+\pi_T\nabla_i\bF^{ij}\wedge\nabla_j\bp+\bF^{ij}\wedge\bh_{ij}\\
&=&-\nabla_i\big(\bF^{ij}\wedge\nabla_j\bp\big)-(\bF^{ij}\cdot\bh^k_i)\nabla_k\bp\wedge\nabla_j\bp\:.
\end{eqnarray*}
The result follows from the symmetry of $\bF^{ij}\cdot\bh_{i}^k$ (cf. Remark~\ref{rem1}). 
$\hfill\blacksquare$\\

\bigskip

We may put together the $\bF$ tensor corresponding to the energy
$$
\mathcal{E}\;:=\;\mathcal{E}_A+\mathcal{E}_0\;=\;\int|\pro d\bH|^2-|\bH\cdot\bh|^2+7|\bH|^4+\mathcal{E}_0\:,
$$
where $\mathcal{E}_0$ is a linear combination of the four conformally invariant energies of order $\mathcal{O}(|\bh_0|^4)$. According to (\ref{GAFA}) and (\ref{tablebraid0}), we have
\be\label{huma1}
\bF\;=\;-\dfrac{1}{2}\Delta_\perp\bH\,g-2(\bH\cdot\bh)\bH-\dfrac{1}{2}\langle\bH\cdot\bh,\bh\rangle g+7|\bH|^2\bH g+\bF_0\:,
\ee
where $\bF_0$ is a normal-valued vector whose trace vanishes. Note that
\be\label{huma2}
\bF^i_i\;=\;-2\Delta_\perp\bH+20|\bH|^2\bH-2\langle\bH\cdot\bh,\bh\rangle\:.
\ee

\bigskip

Among the conformal energies of the type considered in this section, there is of course $\int|W|^2$, with $W$ the Weyl tensor. Recall that 
\be\label{defweyl}
W\;=\;\text{Rm}-\dfrac{1}{2}P\owedge g\:,
\ee
where $P:=\text{Ric}-\frac{1}{6}Rg$ is the Schouten tensor, and $\owedge$ is the Kulkarni-Nomizu product. Using the fact that for a submanifold $\text{Rm}=\dfrac{1}{2}\bh\stackrel{\cdot}{\owedge}\bh$, we obtain after a few elementary manipulations that
$$
W\;=\;\dfrac{1}{2}\bh_0\stackrel{\cdot}{\owedge}\bh_0-\dfrac{1}{2}(h_0^2\owedge g)-\dfrac{1}{6}|\bh_0|^2g\owedge g\:,
$$
where, recall $(h_0^2)^{ij}:=\bh_0^{ik}\cdot(\bh_0)_k^j$ and $\bh_0=\bh-\bH g$. From this expression, one deduces easily that
$$
|W|^2\;=\;2\langle\bh_0^4\rangle-2\text{Tr}_g\bh_0^4-2|h_0^2|^2+\dfrac{1}{3}|\bh_0|^4\:.
$$
We can now read off the expression of $\bF_{|W|^2}$ from the table (\ref{tablebraid0}) and simplify, namely
\bes
\bF_{|W|^2}^{ij}\;=\;8\bh_{ab}W^{aibj}\:.
\ees 
In order to reconstruct $G_{|W|^2}$ from $\bF_{|W|^2}$, we employ the formula given in Remark~\ref{rem1}:
\begin{eqnarray*}
G_{|W|^2}^{ik}&=&-\bh^k_j\cdot\bF^{ij}_{|W|^2}+|W|^2g^{ik}\;\;=\;\;-8W^{aibj}(\bh^k_j\cdot\bh_{ab})+|W|^2g^{ik}\\
&=&-4W^{aibj}R^k_{\:\:ajb}+|W|^2g^{ik}\\
&=&-4W^{iabj}R^k_{\:\:abj}+|W|^2g^{ik}\\
&=&-4W^{iakb}P_{ab}-4W^{iabj}W^k_{\:\:abj}+|W|^2g^{ik}\\
&=&-4W^{iakb}P_{ab}\:,
\end{eqnarray*}
where we have used the Lanczos identity $|W|^{iabj}W^{k}_{\:\:abj}=\frac{1}{4}|W|^2g^{ik}$, which can be found in \cite{Lan}. \\
Altogether, the Noether field associated with the Bach energy is 
\be\label{ojo3}
\bV^i_{|W|^2}\;=\;-4P_{ab}W^{iakb}\nabla_k\bp-8\pro\nabla_j\big(\bh_{ab}W^{aibj}\big)\:.
\ee

\end{document}